\documentclass[10pt,a4paper]{tran-l}

\usepackage {amsmath,amsfonts,amssymb,amsthm,upref,xypic}
\usepackage[poly,curve,arc,dvips]{xy}

\theoremstyle{plain}
\newtheorem{thm}{Theorem}[section]
\newtheorem{prop}[thm]{Proposition}
\newtheorem{cor}[thm]{Corollary}
\newtheorem{lem}[thm]{Lemma}

\theoremstyle{definition}
\newtheorem{defin}[thm]{Definition}
\newtheorem{con}[thm]{Construction}

\theoremstyle{remark}
\newtheorem{rem}[thm]{Remark}

\newtheorem{notation}[thm]{Notation} % This is usually unnumbered

\numberwithin{equation}{section}

\renewcommand{\k}{\mathbf{k}}
\renewcommand{\l}{\mathbf{l}}
\newcommand{\m}{\mathbf{m}}

\renewcommand{\r}{\mathbf{r}}

\newcommand{\sC}{\mathcal{C}}
\newcommand{\sJ}{\mathcal{J}}

\newcommand{\Z}{\mathbb{Z}}
\newcommand{\R}{\mathbb{R}}
\newcommand{\C}{\mathbb{C}}

\newcommand{\N}{\mathbb{N}}

\newcommand{\TTheta}{\mathbf{\Theta}}

\newcommand{\NN}{\widetilde{N}}

\newcommand{\Rbar}{\overline{\mathbb{R}}}

\newcommand{\Rtwoplusbar}{\overline{\R \! \oplus \! \R_+}}

\newcommand{\ubar}{\overline{u}}
\newcommand{\wbar}{\overline{w}}

\renewcommand{\S}{\widetilde{\mathcal{S}}^s}

\newcommand{\Sxk}{S(V) \! \times \! \Delta^k}

\newcommand{\ssigma}{\tilde{\sigma}}

\newcommand{\LS}{\mathbf{LS}}
\newcommand{\LN}{\mathbf{LN}}
\renewcommand{\L}{\mathbf{L}}
\newcommand{\VR}{V \! \oplus \! \R}
\newcommand{\VU}{V \! \oplus \! U}

\newcommand{\Spaces}{\mathrm{Spaces}}
\newcommand{\maps}{\mathrm{maps}}
\newcommand{\hofiber}{\mathrm{hofiber_\ast}}
\newcommand{\hofibery}{\mathrm{hofiber}_y}
\newcommand{\mor}{\mathrm{mor}}
\newcommand{\TOP}{\mathrm{TOP}}
\newcommand{\G}{\mathrm{G}}
\newcommand{\pr}{\mathrm{pr}}
\newcommand{\id}{\mathrm{id}}

\newcommand{\cone}{\mathrm{cone}}
\newcommand{\res}{\mathrm{res}}
\renewcommand{\Bar}{\mathrm{Bar}}
\newcommand{\Stwo}{\widetilde{\mathcal{S}}^{s}_{\Z_2}}

\newcommand{\ra}{\rightarrow}
\newcommand{\co}{\colon\!}

\newcommand{\holimsub}[1]{\begin{array}[t]{cc} \mathrm{holim} \\ [-1mm]
\scriptstyle{#1} \end{array}}

\newcommand{\hocolimsub}[1]{\begin{array}[t]{cc} \mathrm{hocolim} \\
[-1.7mm] \scriptstyle{#1} \end{array}}

\newcommand{\colimsub}[1]{\begin{array}[t]{cc} \mathrm{colim} \\ [-1mm]
\scriptstyle{#1} \end{array}}

\begin{document}

% TOP MATTER

\title[The Block Structure Spaces and Orthogonal Calculus]{The Block Structure
Spaces of Real Projective Spaces and Orthogonal Calculus of
Functors}

\author{Tibor Macko}

\subjclass[2000]{Primary: 57N99, 55P99; Secondary: 57R67}

\keywords{the block structure space, manifold structure, join
construction, orthogonal calculus, the first derivative functor,
splitting problem, surgery}

\address{Mathematisches Institut \\ Universit\"at M\"unster \\
Einsteinstrasse 62 \\ M\"unster, D-48149 \\ Germany \\ and
Matematick\'y \'Ustav SAV \\ \v Stef\'anikova 49 \\ Bratislava,
SK-81473 \\ Slovakia} \email{macko@math.uni-muenster.de}

\begin{abstract}
Given a compact manifold $X$, the set of simple manifold structures
on $X \times \Delta^k$ relative to the boundary can be viewed as the
$k$-th homotopy group of a space $\widetilde{\mathcal{S}}^s (X)$.
This space is called the block structure space of $X$.

We study the block structure spaces of real projective spaces.
Generalizing Wall's join construction we show that there is a
functor from the category of finite dimensional real vector spaces
with inner product to the category of pointed spaces which sends the
vector space $V$ to the block structure space of the projective
space of $V$. We study this functor from the point of view of
orthogonal calculus of functors; we show that it is polynomial of
degree $\leq 1$ in the sense of orthogonal calculus.

This result suggests an attractive description of the block
structure space of the infinite dimensional real projective space
via the Taylor tower of orthogonal calculus. This space is defined
as a colimit of the block structure spaces of projective spaces of
finite dimensional real vector spaces and is closely related to some
automorphisms spaces of real projective spaces.
% \footnote{{\it 2000 Mathematics Subject Classification }57R67}
\end{abstract}

\maketitle
      
\section*{Introduction}
The central objects of study in this paper are the block structure
spaces of certain topological manifolds. These spaces arise as
follows. For a closed topological manifold $X$, a {\it simple
manifold structure} on $X$ is a simple homotopy equivalence $f$ from
some closed manifold $M$ to $X$. Two simple manifold structures $f_0
\colon M_0 \rightarrow X$ and $f_1 \colon M_1 \rightarrow X$ are
{\it equivalent} if there exists a homeomorphism $g \colon M_0
\rightarrow M_1$ and a homotopy from $f_1 \! \circ \! g$ to $f_0$.
The set of equivalence classes is called the {\it simple structure
set} of $X$ and is denoted by $\mathbb{S}^{s}(X)$. More generally,
if $X$ is a compact manifold with boundary, there is a simple
structure set $\mathbb{S}^{s}_\partial(X)$ whose elements are
represented by simple homotopy equivalences of pairs $f \colon (M,
\partial M) \rightarrow (X, \partial X)$ where the restricted map
$\partial M \rightarrow \partial X$ is a homeomorphism.

For a given compact manifold $X$, Quinn constructed in \cite{Quinn}
a space $\S(X)$, the {\it block structure space} of $X$, such that
\begin{equation}
\pi_k(\S(X)) \cong \mathbb{S}^s_{\partial}(X \times \Delta^k).
\end{equation}

This space is related to the automorphism spaces of $X$ via the
homotopy fibration sequence
\begin{equation} \label{two}
\widetilde{\TOP} (X) \rightarrow \G^{s}(X) \rightarrow \S(X).
\end{equation}
Here the space $\widetilde{\TOP}(X)$ is the block version of the
space $\TOP(X)$ of self homeo-morphisms of $X$ and $\G^{s}(X)$ is
the space of simple self homotopy equivalences of $X$. See
\cite{WW}, \cite{BLR}.

The simple structure sets $\mathbb{S}^s(X)$ have been studied in
many special cases. In particular they were calculated for real
projective spaces, lens spaces and complex projective spaces by work
of many authors, especially Wall \cite[chapter 14]{Wall}, Browder
\cite{Br}, Browder-Livesay \cite{BL} and L\'opez de Medrano
\cite{LdM}. A useful tool used by these authors is the {\it join
construction} (see \cite[chapter 14]{Wall}). It provides us with
maps between simple structure sets in different dimensions. For
example in the case of real projective spaces we get maps
$\mathbb{S}^s(\R P^n) \rightarrow \mathbb{S}^s(\R P^{n+m})$. We
generalize the join construction to define maps between the block
structure spaces.

In principle the homotopy type of the block structure space of any
manifold is well understood as the homotopy fiber in the
Sullivan-Wall-Quinn-Ranicki homotopy fibration sequence of surgery
(see \cite[Theorem 2.3.1.]{WW}). Unfortunately, if $m > 2$ this
description of the block structure space does not easily extend to
an illuminating description of the maps $\S(\R P^n) \rightarrow
\S(\R P^{n+m})$ induced by the join construction (Ranicki showed me
how to handle the cases $m = 1,2$). We propose to remedy this defect
by using ideas from functor calculus.

For this purpose we formulate the join construction as a certain
continuous functor. In the case of real projective spaces it is a
functor, which we denote $F$, from the category $\sJ$ of real finite
dimensional vector spaces with inner product to the category
$\Spaces_\ast$ of pointed spaces which sends the vector space $V$ to
the block structure space of $\R P(V)$, the real projective space of
$V$. There are variants of this for the lens spaces and complex
projective spaces. The general methods we use apply to all cases.
However, in the paper we specialize to the case of real projective
spaces. We make remarks on the other cases at the end of the
introduction.

The advantage of the formulation of the join construction as a
continuous functor is that it enables us to use the orthogonal
calculus of functors to study the block structure spaces of real
projective spaces. The orthogonal calculus provides us with a
general theory for studying continuous functors from $\sJ$ (as
above) to $\Spaces_\ast$. It was developed by Weiss in \cite{Weiss}.
He associates to a given continuous functor $E$ the tower of
functors
\[
\cdots \rightarrow T_kE \rightarrow T_{k-1}E \rightarrow \cdots
\rightarrow T_0E,
\]
called the {\it Taylor tower}. For each $k$ the functor $T_kE$ is a
{\it polynomial functor of degree} $\leq k$ in a certain sense and
it is also the best approximation of $E$ by such a functor.
Therefore it should be thought of as the $k$-th Taylor polynomial of
$E$. For example, the functor $T_0E$ is defined by $T_0E(V) =
\mathrm{hocolim}_{n} E(V \oplus\R^n)$. It is the best approximation
of $E$ by an essentially constant functor (all induced maps $T_0E
(V) \rightarrow T_0E(W)$ are homotopy equivalences). Thus the
expansion is at {\it infinity}. The functor $E$ is related to its
Taylor tower via natural transformations $E \rightarrow T_kE$ which
fit together to form a commutative diagram. In favorable
circumstances the maps $E(V) \rightarrow T_kE(V)$ are highly
connected for $k$ large.

The idea of the calculus of functors is to obtain some information
about the functor $E$ by investigating its Taylor tower. For this
purpose we need to estimate the connectivity of the maps $E(V)
\rightarrow T_kE(V)$. The best situation happens if the functor $E$
is polynomial of some degree $\leq k$. Then the maps $E(V)
\rightarrow T_kE(V)$ are homotopy equivalences for all $V \in \sJ$.

A continuous functor $E$ is called polynomial of degree $\leq k$ if
the canonical map
\[
E(V) \rightarrow \holimsub{0 \neq U \subseteq \R^{k+1}} E(\VU)
\]
is a homotopy equivalence for all $V \in \sJ$. The homotopy limit is
over the topological poset of all non-zero vector subspaces of
$\R^{k+1}$ (in section \ref{chap:oc} we give some details about
homotopy limits in this setting). There is a way how to verify that
a given functor is polynomial of degree $\leq k$; in the case $k=1$
we discuss it in detail in the paper.

Once the convergence question is established (that means the
connectivity of maps $E(V) \rightarrow T_kE(V)$ is known), the goal
is to determine the Taylor tower itself. This can be done via an
important theorem in \cite{Weiss} which says that the difference,
$\hofiber(T_kE(V) \rightarrow T_{k-1}E(V))$, between the stages of
the Taylor tower of the functor $E$ is determined by a certain
spectrum $\TTheta E^{(k)}$ with an action of the orthogonal group
$O(k)$. This spectrum should therefore be thought of as the $k$-th
derivative of $E$ at {\it infinity}.

The main theorem of this paper says that the functor $F$ defined
above is essentially polynomial of degree $\leq 1$. It can be
formulated as follows.

\vspace{3mm} \noindent \textbf{Theorem \ref{thm1}.} {\it Let $F
\colon \sJ \rightarrow \Spaces_\ast$ be the functor defined by $V
\mapsto \S(\R P(V))$. Then the canonical maps
\[
F(V) \rightarrow \holimsub{0 \neq U \subseteq \R^2} F(\VU)
\]
are homotopy equivalences if $\dim (V) \geq 6$.} \vspace{3mm}

This implies that the maps $F(V) \rightarrow T_1F(V)$ are homotopy
equivalences if $\dim(V) \geq 6$ (see Corollary \ref{cor1}) and
therefore in high enough dimensions the Taylor tower of $F$ has only
two interesting stages, namely $T_1F$ and $T_0F$. As mentioned
above, the difference between them is measured by the first
derivative of $F$ at infinity, which is a certain spectrum $\TTheta
F^{(1)}$ with an action of the orthogonal group $O(1)$. The homotopy
groups of the spectrum $\TTheta F^{(1)}$ are also determined.

\vspace{3mm} \noindent \textbf{Theorem \ref{thm2}.} {\it We have
\[
\pi_k(\TTheta F^{(1)}) \cong  L_k(1) \mathrm{\;\;for\;\;} k \in \Z,
\]
where the groups $L_k(1)$ are the $4$-periodic $L$-groups from
surgery theory associated to the trivial group.} \vspace{3mm}

There is another spectrum $\L_\bullet$ whose homotopy groups are
$\pi_k (\L_\bullet) \cong L_k(1)$ for $k \in \Z$ (see \cite[chapter
13]{Ran}). The relationship between $\TTheta F^{(1)}$ and
$\L_\bullet$ is not yet clear to us, because we do not know any map
between the two spectra.

The results of Theorems \ref{thm1} and \ref{thm2} can be used to
provide a functorial description of spaces $F(V) \simeq T_1F(V)$ via
the Taylor tower of $F$, if $\dim(V) \geq 6$. For example by
comparing the Taylor tower of $F$ evaluated at $0$ and at $V$ we
obtain the homotopy fibration sequence:
\begin{equation}\label{alt-des}
\Omega^{\infty}[ (S(V)_+ \wedge \TTheta F^{(1)})_{hO(1)} ]
\rightarrow T_1F(0) \rightarrow T_1F(V),
\end{equation}
where $S(V)$ is the unit sphere in $V$ with the antipodal
involution, the subscript $+$ denotes added base point and the
subscript $hO(1)$ denotes the homotopy orbit spectrum. We expect to
determine the homotopy type of $T_1F(0)$ and the first map of the
homotopy fibration sequence (\ref{alt-des}) jointly with Weiss in
\cite{MW}. In any case, this yields a surprising similarity of the
homotopy fibration sequence (\ref{alt-des}) with the
Sullivan-Wall-Quinn-Ranicki homotopy fibration sequence of surgery
theory for $\R P(V)$ (see \cite[Theorem 2.3.1]{WW}). However, an
important difference between the two homotopy fibration sequences is
that (\ref{alt-des}) is functorial in $V$ whereas the surgery
theoretic description is NOT.

The Taylor tower of $F$ can also be used to give a satisfying
description of the block structure spaces of the infinite
dimensional real projective space, that means the colimit
\[
F(\R^\infty) = \colimsub{n \in \N} F(\R^n) \simeq T_0F(0).
\]
Namely, the Taylor tower of $F$ evaluated at $0$ becomes the
homotopy fibration sequence:
\begin{equation} \label{F-infty}
\Omega^{\infty}[ (\TTheta F^{(1)})_{hO(1)} ] \rightarrow T_1F(0)
\rightarrow T_0F(0) \simeq F(\R^\infty).
\end{equation}
If the promised results of \cite{MW} are obtained, this yields the
desired description of $F(\R^\infty)$.

The ultimate goal in the study of the block structure space of a
manifold $X$ is to provide some information about the automorphism
spaces of $X$, the space of self homeomorphisms $\TOP(X)$ and the
space of simple self homotopy equivalences $\G^{s}(X)$. The relation
between these spaces is via the homotopy fibration sequence
(\ref{two}), which involves the block version of the space of self
homeomorphisms $\widetilde{\TOP}(X)$. If $A(X)$ is any of the
automorphism spaces above, the join construction provides us with
maps $A(\R P(\R^n)) \rightarrow A(\R P(\R^m))$ for $n \leq m$ and
the automorphism space $A(\R P^\infty))$ of the infinite dimensional
real projective space is defined as the colimit of $A(\R P(\R^n))$
over $n \in \N$. The join maps commute with the maps in (\ref{two})
and so we obtain the homotopy fibration sequence
\begin{equation} \label{two-stable}
\widetilde{\TOP} (\R P^\infty)) \rightarrow \G^{s}(\R P^\infty))
\rightarrow \S(\R P^\infty)) = F(\R^\infty).
\end{equation}
Closely related to these spaces are also certain equivariant
automorphism spaces discussed in a more detail below. We remark that
the space $\mathrm{colim}_{n \in \N} \G^s_{\Z_2}(S(n \cdot V))$
defined below is a trivial double cover of $\G^s(\R P^\infty))$ and
in \cite{BS} it is shown that this space is homotopy equivalent to
$Q(\R P^\infty_+)$, where $Q = \Omega^\infty \Sigma^\infty$ and the
subscript $+$ denotes the added base point. Therefore any
information about the homotopy type of $F(\R^\infty)$ gives us also
information about $\widetilde{\TOP}(\R P^\infty))$ which is closely
related to $\TOP(\R P^\infty))$ (see \cite{WW} for more information
about the relation between various automorphism spaces, also note
that for real projective spaces the word simple can be dropped as
$\mathrm{Wh} (\Z_2)= \{0\}$).

There is also a more direct way to relate the space $F(\R^\infty)$
to certain honest (not block) equivariant automorphism spaces as
follows. Let $W$ be the regular representation of the group $\Z_2$.
Denote by $\G_{\Z_2}(S(n \cdot W))$ the topological monoid of
equivariant self homotopy equivalences of the unit sphere $S(n \cdot
W)$ of the induced representation $n \cdot W = \R^n \otimes W$ and
by $\TOP_{\Z_2}(S(n \cdot W))$ the topological group of equivariant
self homeomorphisms of $S(n \cdot W)$ . The join map induces
inclusions $\G_{\Z_2}(S(n \cdot W)) \hookrightarrow \G_{\Z_2}(S(m
\cdot W))$ and $\TOP_{\Z_2}(S(n \cdot W)) \hookrightarrow
\TOP_{\Z_2}(S(m \cdot W))$ for $n \leq m$ and there are also
inclusions $\TOP_{\Z_2}(S(n \cdot W)) \hookrightarrow \G_{\Z_2}(S(n
\cdot W))$ which commute with the join map. Therefore we get the
maps between quotients induced by the join construction. After
passing to the colimit over $n \in  \N$ there is the following
splitting
\begin{equation} \label{equiv-spaces}
\colimsub{n \in \N} \G_{\Z_2}(S(n \cdot W))/\TOP_{\Z_2}(S(n \cdot
W)) \simeq F(\R^\infty)_1 \times \G/\TOP,
\end{equation}
which relates the left hand side to the space $F(\R^\infty)$. Here
$F(\R^\infty)_1$ is the base point component. For the classical
space $\G/\TOP$ see e.g. \cite[section 1]{WW}.

Here is a sketch of how the splitting (\ref{equiv-spaces}) is
obtained. Let $W=V \oplus U$ where $V$ is the free part of the
representation and $U$ is the trivial part. There is the splitting
\begin{equation} \label{htpy-splitting}
\colimsub{n \in \N} \G_{\Z_2}(S(n \cdot W)) \simeq \colimsub{n \in
\N} \G_{\Z_2}(S(n \cdot V)) \times \colimsub{n \in \N} \G (S(n \cdot
U)).
\end{equation}
See \cite[Proposition 3.5.]{Crabb} for the statement on the level of
homotopy groups. Similarly there is the splitting
\begin{equation} \label{top-stable-splitting}
\colimsub{n \in \N} \TOP_{\Z_2}(S(n \cdot W)) \simeq \colimsub{n \in
\N} \widetilde{\TOP}_{\Z_2}(S(n \cdot V)) \times \colimsub{n \in \N}
\TOP (S(n \cdot U)).
\end{equation}
This is obtained by passing to the colimit over $n$ with the
splittings
\begin{equation} \label{top-splitting}
\TOP_{\Z_2}(S(n \cdot W)) \simeq \TOP_{\Z_2}(S(n \cdot W), \;
\mathrm{rel} \; S(n \cdot U)) \times \TOP (S(n \cdot U))
\end{equation}
where $\; \mathrm{rel} \; S(n \cdot U)$ means that the equivariant
self homeomorphisms of $S(n \cdot W)$ are fixed to be the identity
on $S(n \cdot U) \hookrightarrow S(n \cdot W)$. The splitting
(\ref{top-splitting}) is analogous to the splitting in \cite[Lemma
3.4.]{Crabb}. The second factor in (\ref{top-stable-splitting}) is
the colimit of the second factors of (\ref{top-splitting}). The
first factor of (\ref{top-stable-splitting}) is obtained by relating
the first factor of (\ref{top-splitting}) to the equivariant bounded
self homeomorphisms of $S(n \cdot V) \times (n \cdot U)$ using a
trick of \cite[section 3]{AH} and then relating the bounded self
homeomorphisms to block self homeomorphisms as in \cite[Proposition
1.4.2.]{WW}. Finally, these equivariant automorphism spaces of $S(n
\cdot V)$ map into automorphism spaces of $\R P(n \cdot V)$ by
passing to quotients, and this map is a trivial double cover. The
splitting (\ref{equiv-spaces}) then follows from (\ref{two-stable})
because the join maps commute with the splittings
(\ref{htpy-splitting}), (\ref{top-stable-splitting}) and we have
\begin{eqnarray*}
\G_{\Z_2} (S(n \cdot V))/\widetilde{\TOP}_{\Z_2} (S(n \cdot V))
 & \cong & \\ \G (\R P(n \cdot V))/\widetilde{\TOP} (\R P(n \cdot
V)) & \cong & \S_1 (\R P (n \cdot V)).
\end{eqnarray*}

%%%%%%%%%%%%%%%%%%%%%%%%%%%%%%%%%%%%%%%%%%%%%%%%%%%%%%%%%
%CORRECTION%
%%%%%%%%%%%%%%%%%%%%%%%%%%%%%%%%%%%%%%%%%%%%%%%%%%%%%%%%%

Now we would like to briefly
discuss the study of the block structure spaces of lens spaces and
complex projective spaces. Using the join construction we obtain in
these cases a functor from the category of complex vector spaces
with inner product to $\Spaces_\ast$ which sends the vector space
$V$ to $\S(L(V))$ or $\S(\C P(V))$. Here we can employ the unitary
calculus which is the complex analogue of orthogonal calculus, and
we can try to run the study parallel to the case of real projective
spaces. The case of lens spaces was studied by Madsen and Rothenberg
in \cite{MR}, where they show that the join maps induce isomorphisms
on homotopy groups (they also determine the homotopy groups, and
later Madsen in \cite{Mad} determines the homotopy type localized at
$2$ and away from $2$). Hence, from our point of view, the
corresponding functor is polynomial of degree $0$. The case of
complex projective spaces is not covered in \cite{MR}, so the
parallel study would be interesting. We note that the main
ingredient in the proof of Theorems \ref{thm1} and \ref{thm2} is the
codimension $1$ surgery theory, whereas in the case of complex
projective spaces we have to use codimension $2$ surgery theory,
which is considerably more complicated (see \cite[chapter
7.8]{Ran1}). At the moment we are unable to carry over the whole
proofs of Theorems \ref{thm1}, \ref{thm2}. We hope to do this in the
future.

The paper is organized as follows. In section  \ref{blockspace} we
introduce the block structure spaces. In section
\ref{chap:suspension} we construct the functor $F$. In section
\ref{chap:oc} we give an overview of the basic constructions from
orthogonal calculus that we use in the paper. In particular for a
continuous functor $E \colon \sJ \rightarrow \Spaces_\ast$ we recall
a necessary and sufficient condition for $E$ to be polynomial of
degree $\leq 1$ (Theorem \ref{poly1}). This involves studying a
certain functor $E^{(1)} \colon \sJ \rightarrow \Spaces_\ast$ with
some additional structure, which is called the first derivative
functor of $E$. Section \ref{chap:F1} contains background material
about codimension $1$ surgery theory, which is used to give a
description of the first derivative functor $F^{(1)}$. Up to now we
claim no originality. The description of $F^{(1)}$ will then be the
main tool in the proofs of Theorems \ref{thm1} and \ref{thm2} which
are given in section \ref{chap:ThetaF1}. The proofs are our
contribution, and they are the core of the paper. We note that
section \ref{chap:oc} is a short survey of orthogonal calculus,
however, for complete information we refer the reader to the
original source \cite{Weiss}. Throughout the paper we also assume
basic knowledge of surgery theory (that means the book \cite{Wall}
up to chapter 10), but we provide some background about the less
standard parts, such as codimension 1 surgery theory.

\begin{notation}
For $V \in \sJ$ by $\R P(V)$ we denote the real projective space of
$V$, by $S(V,r)$ we denote the sphere with radius $r \geq 0$ in $V$;
and $S(V) = S(V,1)$.
\end{notation}

%%%%%%%%%%%%%%%%%%%%%%%%%%%%%%%%%%%%%%%%%%%%%%%%%%%%%%%%%%%%%%
%End of CORRECTION%
%%%%%%%%%%%%%%%%%%%%%%%%%%%%%%%%%%%%%%%%%%%%%%%%%%%%%%%%%%%%%%

%%%%%%%%%%%%%%%%%%%%%%%%%%%%%%%%%%%%%%%%%%%%%%%%%%%%%%%%%%%%%%

\section{The Block Structure Space of a Manifold} \label{blockspace}

%%%%%%%%%%%%%%%%%%%%%%%%%%%%%%%%%%%%%%%%%%%%%%%%%%%%%%%%%%%%%%
%CORRECTION%
%%%%%%%%%%%%%%%%%%%%%%%%%%%%%%%%%%%%%%%%%%%%%%%%%%%%%%%%%%%%%%

In this section we introduce the
block structure space of a closed manifold $X$ in general and we
also introduce a method for constructing other models for the block
structure spaces. The method will be used in section
\ref{chap:suspension} to obtain a functorial model of the block
structure spaces of the real projective spaces.

%%%%%%%%%%%%%%%%%%%%%%%%%%%%%%%%%%%%%%%%%%%%%%%%%%%%%%%%%%%%%%
%End of CORRECTION%
%%%%%%%%%%%%%%%%%%%%%%%%%%%%%%%%%%%%%%%%%%%%%%%%%%%%%%%%%%%%%%

In the introduction we have defined for a given closed manifold $X$
its simple structure set $\mathbb{S}^s(X)$ and more generally for a
compact manifold with boundary $(X,\partial X)$ the simple structure
set $\mathbb{S}^s_\partial(X)$. These are the basic objects of study
of surgery theory which have been calculated in many special cases
(see for example \cite{Wall}, \cite{LdM}).

For a compact manifold $X$ the structure sets
$\mathbb{S}^s_\partial(X \times \Delta^k)$ become groups when $k
\geq 1$. The block structure space of $X$ is designed to be a space
such that
\begin{equation}
\pi_k(\S(X)) \cong \mathbb{S}^{s}_{\partial}(X \times \Delta^k).
\label{htpygroups}
\end{equation}
It is constructed as a $\Delta$-set (see \cite{RS}). If $X$ is a
closed manifold then a $k$-simplex of $\S(X)$ is a simple homotopy
equivalence $f \colon M \rightarrow X \times \Delta^k$ between
$(n+k)$-dimensional manifold $(k+2)$-ads, where $n = \dim(X)$ (the
notion of a manifold $(k+2)$-ad is discussed below). The face
operators are given by restrictions and the base point in dimension
$k$ is $\id \colon X \times \Delta^k \rightarrow X \times \Delta^k$.

The construction of the space $\S(X)$ was given by Quinn in
\cite{Quinn}. As indicated in the introduction the main motivation
for studying this space comes from the fact that it is closely
related to the automorphism spaces of $X$. We refer the reader to
the survey article \cite{WW} for more information about this
application. In \cite[Theorem 2.3.1.]{WW} it is also stated that
$\S(X)$ fits into the homotopy fibration sequence of surgery theory
of Sullivan-Wall-Quinn-Ranicki, which in principle determines its
homotopy type.

%%%%%%%%%%%%%%%%%%%%%%%%%%%%%%%%%%%%%%%%%%%%%%%%%%%%%%%%%%%%%%
%CORRECTION%
%%%%%%%%%%%%%%%%%%%%%%%%%%%%%%%%%%%%%%%%%%%%%%%%%%%%%%%%%%%%%%

As already mentioned, when $X$ is
the real projective space $\R P(V)$ for some $V \in \sJ$, we obtain
an additional structure of maps between the block structure spaces
of real projective spaces in different dimensions. These maps,
called the {\it join maps}, were described on the homotopy groups by
Wall \cite{Wall} for $\pi_0$ and by Madsen and Rothenberg \cite{MR}
for $\pi_k$, when $k \geq 1$. We generalize these maps to the
level of spaces in a functorial way. For this we need another model
of the block structure space of $\R P(V)$, that means a space
homotopy equivalent to $\S (\R P(V))$. The method for constructing
such a space, which we describe in this section, is based on a
careful analysis of Quinn's model. The desired model is then 
constructed in two steps. In this section we construct an
intermediate model $\S_{\Z_2} (S(V))$, which also admits the join
maps, as explained in \cite{MR}. However, these do not define a
functor because they do not behave well w.r.t. the composition.
Then, in the next section, we construct the functorial model
$\S_{\Z_2}(V)$.

\begin{rem} \label{models}
Each of the models $\S (\R P(V))$, $\Stwo (S(V))$, $\Stwo (V)$ has
convenient properties in different situations. As they are all
homotopy equivalent (at least in high dimensions) and as many of the
results we need are only `up to homotopy' we switch between these models throughout the paper according to our needs.
\end{rem}

%%%%%%%%%%%%%%%%%%%%%%%%%%%%%%%%%%%%%%%%%%%%%%%%%%%%%%%%%%%%%%
%End of CORRECTION%
%%%%%%%%%%%%%%%%%%%%%%%%%%%%%%%%%%%%%%%%%%%%%%%%%%%%%%%%%%%%%%

Firstly we elaborate on the definition of a manifold $(k+2)$-ad. It
is a generalization of a cobordism between manifolds, namely a
cobordism is a special manifold $3$-ad. According to \cite{Mil} (or
\cite{Luck}) a cobordism between two $n$-dimensional manifolds
$M_0$, $M_1$ is a $5$-tuple $(M,M_0,M_1,h_0,h_1)$, where $M$ is an
$(n+1)$-dimensional manifold with two disjoint pieces of boundary
$\partial_0 M$, $\partial_1 M$ and $h_i \colon M_i \rightarrow
\partial_i M$ for $i=0,1$ are homeomorphisms of $n$-dimensional
manifolds. A manifold $(k+2)$-ad is a generalization where the
manifold $M$ has $(k+1)$ pieces of boundary which fit together in a
nice way. We formalize this notion below.

More precisely we give an abstract construction of a $\Delta$-set
which consists of objects formally resembling manifold $(k+2)$-ads
as an associated $\Delta$-set to a certain $\Delta$-groupoid. This
is done in Construction \ref{abstract} and then applied to the
special case of the block structure spaces in Construction
\ref{blockdef}.

We recall basic definitions for $\Delta$-sets of \cite{RS} (see also
\cite[chapter 11]{Ran}). A $\Delta$-set $K$ is a functor $K \colon
\Delta^{op} \rightarrow \mathrm{Sets}$, where $\Delta$ is the
category with objects finite totally ordered sets $\k = \{ 0,
\ldots, k \}$ and morphisms injective monotone functions $\alpha
\colon \m \rightarrow \k$ for $m \leq k$. Thus a $\Delta$-set $K$
can be described as a collection of sets $\{K_k\}_{k \in \N}$
together with face maps $\partial_\alpha \colon K_k \rightarrow
K_m$. The standard $k$-simplex is a $\Delta$-set in the following
way: the set of $m$-simplices $\Delta^k_m$ is the set of monotone
injective functions $\sigma \colon \m \rightarrow \k$ and for $\beta
\colon \l \rightarrow \m$ the face map $\partial_\beta \colon
\Delta^k_m \rightarrow \Delta^k_l$ is defined by
$\partial_\beta(\sigma) = \sigma \! \circ \! \beta$. Also $\alpha
\colon \m \rightarrow \k$ defines a map of $\Delta$-sets $\Delta^m
\rightarrow \Delta^k$. More generally, a $\Delta$-groupoid $K$ is a
functor $K \colon \Delta^{op} \rightarrow \mathrm{Groupoids}$, that
means for each $k$ the set of $k$-simplices $K_k$ is a groupoid (a
category where all morphisms are isomorphisms) and the face maps
$\partial_\alpha$ are functors. Note that the standard simplex
$\Delta^k$ is also a $\Delta$-groupoid in a trivial way. A
$\Delta$-set $K$ is {\it pointed} if for every $k \in \N$ there is a
base simplex $\ast \in K_k$ such that $\partial_\alpha \ast = \ast$
for all $\alpha$.

Now we recall the following well known construction for the special
case of a $\Delta$-groupoid.

\begin{con} \cite{Dwyer}
Let $K$ be a $\Delta$-groupoid. The {\it Grothendieck construction}
$Gr(K)$ of $K$ is the category with objects pairs $(\k,x)$, where
$\k \in \Delta$, $x \in K_k$, and morphisms $(\k,x) \rightarrow
(\m,y)$ are pairs $(\alpha, h)$, where $\alpha \colon \m \rightarrow
\k$ is a morphism in $\Delta$ and $h \colon
\partial_\alpha x \rightarrow y$ is a morphism in $K_m$. We also
define the degree of the object $(\k,x)$ to be $k \in \N$. Note that
the construction is functorial. In particular for $\alpha \colon \m
\rightarrow \k$ we have a functor $Gr(\alpha) \colon Gr(\Delta^m)
\rightarrow Gr(\Delta^k)$.
\end{con}

\begin{con} \label{abstract}
Let $K$ be a $\Delta$-groupoid. We define the {\it associated
$\Delta$-set} of $K$ to be the $\Delta$-set $K^\sharp$ defined as
follows. The set of $k$-simplices of the $\Delta$-set $K^\sharp$ is
the set of degree preserving functors $\sigma \colon Gr(\Delta^k)
\rightarrow Gr(K)$ which satisfy the extra condition that for $\beta
\colon \l \rightarrow \m$ the morphism $\sigma (\beta, id)$ equals
$(\beta, h_\beta)$ for some $h_\beta \in K_l$. For $\alpha \colon \m
\rightarrow \k$ the face map $\partial^\sharp_\alpha \colon
K^\sharp_k \rightarrow K^\sharp_m$ is defined by
$\partial^\sharp_\alpha \sigma = \sigma \! \circ \! Gr(\alpha)$.
\end{con}

A $k$-simplex $\sigma \in K^\sharp$ can be described by a list
$\sigma = (x_\alpha, h_\alpha)$, where $x_\alpha \in K_m$ and
$h_\alpha \colon \partial_\alpha x \rightarrow x_\alpha$ is an
isomorphism in $K_m$ for $\alpha \colon \m \rightarrow \k$ a
morphism in $\Delta$, satisfying various compatibility conditions.
Here $x=x_\id$. Note that if $K$ is a trivial $\Delta$-groupoid (all
morphisms are the identity), then a $k$-simplex of $K^\sharp$ is
just a collection $\sigma = (x_\alpha, \id)$, where $x_\alpha =
\partial_\alpha x$. Thus we recover the underlying $\Delta$-set
of $K$. If $K$ is a non-trivial $\Delta$-groupoid a $k$-simplex of
$K^\sharp$ is a collection $\sigma = (x_\alpha, h_\alpha)$, where
objects $x_\alpha$ are just isomorphic to $\partial_\alpha x$.

\begin{rem} \label{rem:sharp}
The assignment $K \mapsto K^\sharp$ is a functor from
$\Delta$-groupoids to $\Delta$-sets. Suppose that $K$ is a
topological $\Delta$-groupoid, that means the category $K_k$ is a
topological category for each $k$ and the face operators are
continuous functors. Then $K^\sharp_k$ can also be endowed with a
topology (the subspace topology of the product topology) and
$K^\sharp$ becomes a $\Delta$-space, also in a functorial way. The
geometric realization then converts the $\Delta$-set (or the
$\Delta$-space) $K^\sharp$ into a genuine space again in a
functorial way. In the notation we do not distinguish between the
underlying $\Delta$-set (or $\Delta$-space) and its geometric
realization. We note that all these constructions have pointed
versions.
\end{rem}

%%%%%%%%%%%%%%%%%%%%%%%%%%%%%%%%%%%%%%%%%%%%%%%%%%%%%%%%%%%%%%
%CORRECTION%
%%%%%%%%%%%%%%%%%%%%%%%%%%%%%%%%%%%%%%%%%%%%%%%%%%%%%%%%%%%%%%

\begin{defin}
A $k$-{\it block} $M$ is a manifold with boundary such that the
boundary $\partial M$ is decomposed into codimension $0$
submanifolds $\partial_i M$, for $i \in \k$, whose intersections
$\partial_i M \cap \partial_j M$ are codimension $0$ submanifolds of
the boundaries of $\partial_i M$ and $\partial_j M$, and so on. For
$\alpha \colon \m \rightarrow \k$ a morphism in $\Delta$ we denote
$\partial_\alpha M = \bigcap_{i \in \mathrm{im}(\alpha)} \partial_i
M$. A {\it block map} $f \colon M \rightarrow N$ between two
$k$-blocks is a map, such that $f (M_\alpha) \subset N_\alpha$ for
$\alpha \in \k$. A $k$-block is {\it special} if $\bigcap_{i \in \k}
\partial_i M = \emptyset$.
\end{defin}

An example of a special $k$-block is $X \times \Delta^k$ for $X$ a
closed manifold. Now we are ready for the definition of the block
structure space $\S(X)$.

\begin{con} \label{blockdef}
Let $X$ be a closed $n$-dimensional manifold. The block structure
space $\S(X)$ of $X$ is the associated pointed $\Delta$-set of the
pointed $\Delta$-groupoid whose $k$-simplices are pairs $(M,f)$, where $M$ is a special $k$-block, and a block map $f \colon M \rightarrow X \times \Delta^k$
is a simple homotopy equivalence between $(n+k)$-dimensional
$k$-blocks. Face operators $\partial_\alpha$ are obtained by taking
the face $\partial_\alpha M$ and restricting $f$ to
$f|_{\partial_\alpha M} \colon \partial_\alpha M \rightarrow X
\times \partial_\alpha \Delta^k$. An isomorphism between $(M,f)$ and
$(M',f')$ is a block homeomorphism $h \colon M \rightarrow M'$ such that $f = f' \circ h$. The base point in dimension
$k$ is $\id \colon X \times \Delta^k \rightarrow X \times \Delta^k$.
\end{con}

%%%%%%%%%%%%%%%%%%%%%%%%%%%%%%%%%%%%%%%%%%%%%%%%%%%%%%%%%%%%%%
%End of CORRECTION%
%%%%%%%%%%%%%%%%%%%%%%%%%%%%%%%%%%%%%%%%%%%%%%%%%%%%%%%%%%%%%%

A $k$-simplex of $\S(X)$ can be described by a list $\,((M_\alpha,
f_\alpha), h_\alpha)$, where for $\alpha \colon \m \rightarrow \k$
the $M_\alpha$ is a special $m$-block, the map $f_\alpha \colon
M_\alpha \rightarrow X \times \partial_\alpha \Delta^k$ is a simple
homotopy equivalence of $m$-blocks, $h_\alpha \colon M_\alpha
\rightarrow \partial_\alpha M$ is a homeomorphism of $m$-blocks and
the maps $f_\alpha$, $h_\alpha$ satisfy various compatibility
conditions. The list $(M_\alpha, h_\alpha)$ is the correct
formalization of a manifold $(k+2)$-ad. In the case $k=1$, the
manifold $3$-ad is a cobordism in the sense of \cite{Mil},
\cite{Luck} as described earlier.

The block structure space $\S(X)$ as defined above is a Kan
$\Delta$-set for any manifold $X$. For Kan $\Delta$-sets there is a
convenient combinatorial formula for homotopy groups (see
\cite[chapter 11]{Ran}) and using it we immediately obtain the
isomorphisms (\ref{htpygroups}).

%%%%%%%%%%%%%%%%%%%%%%%%%%%%%%%%%%%%%%%%%%%%%%%%%%%%%%%%%%%%%%
%CORRECTION%
%%%%%%%%%%%%%%%%%%%%%%%%%%%%%%%%%%%%%%%%%%%%%%%%%%%%%%%%%%%%%%

The following lemma provides us
with a method for constructing other models of the block structure
space.

\begin{lem} \label{skelet}
Let $L \rightarrow K$ be a map of pointed $\Delta$-groupoids such
that $L^\sharp$ and $K^\sharp$ are Kan $\Delta$-sets. If $L_k
\rightarrow K_k$ induces a bijection on equivalence classes of
objects for every $k$, then the induced map $L^\sharp \rightarrow
K^\sharp$ is a homotopy equivalence.
\end{lem}

\begin{proof}
Using the combinatorial expression of $\pi_k$ it is easy to see that
the map $L^\sharp \rightarrow K^\sharp$ induces isomorphisms on
homotopy groups.
\end{proof}

Now we specialize to the real projective spaces. Note first that any
homotopy equivalence $f \colon M \rightarrow \R P(V) \times
\Delta^k$ is simple because $\mathrm{Wh}(\Z_2) = \{0\}$ (see
\cite{Luck}). Therefore we drop the word simple when dealing with
the real projective spaces. Here is the promised intermediate model
$\Stwo (S(V))$ for the block structure space of $\R P(V)$.

\begin{con}\label{F-spheres}
\label{doublecovers} Let $\Stwo(S(V))$ be the associated pointed
$\Delta$-set of the pointed $\Delta$-groupoid whose $k$-simplex is a
pair $(T,f)$, where:
\begin{itemize}
 \item $T$ is a free involution on $\Sxk$, which, as a map, is a block map,
 \item $f \colon \Sxk \rightarrow \Sxk$ is a block map which is an
 equivariant homotopy equivalence w.r.t. $T$ on the source and $T_a$
 on the target, where $T_a$ is the product of to the antipodal
 involution on $V$ (restricted to $S(V)$) with the identity on $\Delta^k$.
\end{itemize}
Face maps are defined by restriction and an isomorphism between
$(T,f)$ and $(T',f')$ is an equivariant block homeomorphisms $h
\colon \Sxk \rightarrow \Sxk$ such that $f = f' \circ h$. The base
point is the pair $(T_a,\id)$. A pair $(T,f)$ which represents a
$k$-simplex in $\Stwo (S(V))$ will be shortly referred to as a {\it
structure} on $S(V) \times \Delta^k$
\end{con}

There is a map $\Stwo(S(V)) \rightarrow \S(\R P(V))$ defined on the
underlying $\Delta$-groupoids by taking quotients.

\begin{prop} \label{modelone}
The map $\Stwo(S(V)) \rightarrow \S ( \R P(V))$ is a homotopy
equivalence if $\dim (V) \geq 5$.
\end{prop}

\begin{proof}[Proof of Proposition \ref{modelone}]
Using Lemma \ref{skelet} it is enough to show that for each $k \in
\N$ the map from the $k$-simplices of the underlying groupoid of
$\Stwo(S(V))$ to the $k$-simplices of the underlying groupoid of
$\S(\R P(V))$ induces an isomorphism on the equivalence classes of
objects. We proceed by induction on $k$. From the generalized
Poincar\'e conjecture in $\dim (V) \geq 5$ we have that any manifold
$M$ homotopy equivalent to $\R P(V)$ has the universal cover
homeomorphic to $S(V)$ which proves the statement for $k = 0$. For
the induction step we apply the $h$-cobordism theorem which implies
that the universal cover of any manifold $(k+2)$-ad $M$ mapping via
a homotopy equivalence to $\R P(V) \times \Delta^k$ is homeomorphic
to $S(V) \times \Delta^k$.
\end{proof}

%%%%%%%%%%%%%%%%%%%%%%%%%%%%%%%%%%%%%%%%%%%%%%%%%%%%%%%%%%%%%%
%End of CORRECTION%
%%%%%%%%%%%%%%%%%%%%%%%%%%%%%%%%%%%%%%%%%%%%%%%%%%%%%%%%%%%%%%

%%%%%%%%%%%%%%%%%%%%%%%%%%%%%%%%%%%%%%%%%%%%%%%%%%%%%%%%%%%%%%%%%%

\section{The Join Construction}\label{chap:suspension}

%%%%%%%%%%%%%%%%%%%%%%%%%%%%%%%%%%%%%%%%%%%%%%%%%%%%%%%%%%%%%%
%CORRECTION%
%%%%%%%%%%%%%%%%%%%%%%%%%%%%%%%%%%%%%%%%%%%%%%%%%%%%%%%%%%%%%%

In this section we construct a
continuous functor $F \colon \sJ \rightarrow \Spaces_\ast$, such
that for an object $V \in \sJ$ with $\dim(V) \geq 5$, the
space $F(V)$ is homotopy equivalent to the space $\S (\R P(V))$ and
for a morphism $\xi \in \sJ$ the map $F(\xi)$ is a generalization of
the join construction of Wall.

\begin{con}\label{F-functor}
Given an object $V \in \sJ$, let $\Stwo (V)$ be the associated
pointed $\Delta$-set of the pointed $\Delta$-groupoid whose
$k$-simplex is a pair $(T,f)$ where:
\begin{itemize}
\item $T$ is an involution on $V \times \Delta^k$, which, as a map, is
a block map, it leaves $S(V,r) \times \Delta^k$ invariant for each
$r \geq 0$, its restriction to $S(V,r) \times \Delta^k$ is free for
$r > 0$, and its restriction to $0 \times \Delta^k$ is trivial,
\item $f \colon V \times \Delta^k \rightarrow V \times
\Delta^k$ is a block map, which is equivariant w.r.t. $T$ on the
source and $T_a$ on the target, it leaves $S(V,r) \times \Delta^k$
invariant for each $r \geq 0$, the restriction of $f$ to $S(V,r)
\times \Delta^k$ is an equivariant homotopy equivalence for $r > 0$,
and the restriction of $f$ to $0 \times\Delta^k$ is a homeomorphism.
\end{itemize}
Face maps are defined by restriction, an isomorphism from $(T,f)$ to
$(T',f')$ is an equivariant block homeomorphism $h \colon V \times
\Delta^k \rightarrow V \times \Delta^k$ such that $f = f' \circ h$,
and the base point is the pair $(T_a,\id)$. A pair $(T,f)$ which
represents a $k$-simplex in $\Stwo (V)$ will be shortly referred to
as a {\it structure} on $V \times \Delta^k$.
\end{con}

\begin{con}
Given a morphism $\xi \co V \ra W$ in $\sJ$, let $\Stwo (\xi) \colon
\Stwo (V) \rightarrow \Stwo (W)$ be a map defined as follows. The
linear inclusion $\xi$ induces a direct sum decomposition $W =
\xi(V) \oplus \xi(V)^\perp$ and the product decomposition $W \times
\Delta^k = \xi(V) \times \Delta^k \times \xi(V)^\perp$. First we
obtain the structure $(\xi \circ T \circ \xi^{-1}, \xi \circ f \circ
\xi^{-1})$ on $\xi(V) \times \Delta^k$. The desired structure on $W
\times \Delta^k$ is defined by
\[
(\xi \circ T \circ \xi^{-1}, \xi \circ f \circ \xi^{-1}) \times
(T_a,\id)_{\xi(V)^\perp}.
\]
\end{con}

\begin{thm}
Given an object $V \in \sJ$, let $F(V)$ be the geometric realization
of $\Stwo(V)$, and given a morphism $\xi \in \sJ$, let $F(\xi)$ be the geometric realization of $\Stwo(\xi)$. Then $F$ defines a functor from $\sJ$ to $\Spaces_\ast$.
\end{thm}

The statement is easy to verify and we leave it for the reader. In
order to show that the space $\Stwo (V)$ defines another model for
the block structure space of $\R P(V)$ we use the method of the
previous section to prove the following proposition.

\begin{prop}\label{spheres-vs-vector-spaces}
For each $V \in \sJ$ there exist maps
\[
\res \co \Stwo (V) \ra \Stwo (S(V)), \quad \text{and} \quad \cone
\co \Stwo (S(V)) \ra \Stwo (V)
\]
which are homotopy inverses of each other.
\end{prop}

\begin{rem}
As already indicated, throughout the paper we also work with the
models $\S(\R P(V))$ and $\Stwo (S(V))$, which is legal for the
statements which are `up to homotopy'. Similarly, we work with the
join maps on these models. For example in section \ref{chap:ThetaF1}
the model for the join map is $(\res \circ F(\xi) \circ \cone)$.
\end{rem}

\begin{con}
The maps $\res$ and $\cone$ are induced by the maps between the
corresponding pointed $\Delta$-groupoids. Since we think of $S(V)
\times \Delta^k = S(V,1) \times \Delta^k$ as of a subspace of $V
\times \Delta^k$, we can define
\[
\res \co \Stwo (V) \ra \Stwo (S(V))
\]
by $(T,f) \mapsto (T|_{(S(V) \times \Delta^k)},f|_{(S(V) \times
\Delta^k)})$.

The map
\[
\cone \co \Stwo (S(V)) \ra \Stwo (V)
\]
is defined by a certain `iterated cone construction' which we now
describe. Let $(T,f)$ be a structure on $S(V) \times \Delta^k$ and let
$(T',f') = \cone(T,f)$, which is a structure on $V \times \Delta^k$.
We describe a general procedure for constructing a block self-map of
$V \times \Delta^k$, which preserves $S(V,r) \times \Delta^k$ for
each $r \geq 0$, from a block self-map of $S(V) \times \Delta^k$. The
maps $T'$, $f'$ are constructed from $T$, $f$ using this procedure.

\medskip\noindent\textbf{The case $k=0$.} Let $X$ be a topological space. The {\it open
cone} on $X$ is the space $C^+X$ whose points are linear expressions
$rx$, where $r \in [0,\infty)$. This includes the identifications
$0=0x=0y$. The point $0$ is called the {\it cone point}. We have a
decomposition $C^+X = CX \cup C^{1+}X$, with $CX = \{ rx \; | \; r
\in [0,1]\}$, and $C^{1+}X = \{ rx \; | \; r \in [1,\infty)\}$.  For
$r \geq 0$ we distinguish the subspaces $r \cdot X =\{rx \; | \; x
\in X\}$ and we identify $X = 1 \cdot X$. Sometimes we want to give
the cone point some other name than $0$, e.g. $c$. In that case the
points of the subspace $CX$ are linear combinations $rx+(1-r)c$ for
$x \in X$ and $r \in [0,1]$. Note that, if $X$ is $S(V)$, then $C^+ X$ can be identified with $V$, so that $r \cdot X$ is identified with $S(V,r)$ for $r \geq 0$.

Let $h \colon X \rightarrow Y$ be a map. Then the rule $rx \mapsto
rh(x)$ defines an extension $C^+(h)$ of $h$ such that $C^+(h) (r
\cdot X) \subseteq r \cdot Y$.

\medskip\noindent\textbf{The cases $k>0$.} Let $h \colon X \times
\Delta^k \rightarrow Y \times \Delta^k$ be a block map. We define
its extension, a block map $C^+ (h) \colon C^+X \times \Delta^k
\rightarrow Y \times \Delta^k$, such that $C^+ (h) (r \cdot X \times
\Delta^k) \subseteq r \cdot Y \times \Delta^k$.

We have the decomposition $C^+ X \times \Delta^k = CX \times
\Delta^k \cup C^{1+}X \times \Delta^k$. On $C^{1+}X \times \Delta^k$
the extension is defined as in the case $k=0$, i.e. if $h(x,s) =
(h(x,s)_1,h(x,s)_2)$, then $C^+ (h) \colon (rx,s) \mapsto
(rh(x,s)_1,h(x,s)_2)$.  On $CX \times \Delta^k$ we use the following
`iterated cone construction'.

Starting with $X \times \Delta^k$ we define the {\it iterated cone}
$C_k (X \times \Delta^k)$ by
\begin{equation}\label{iterated-cone}
C_k (X \times \Delta^k) = C \big( X \times \Delta^k \cup \bigcup_{\tau} C_{k-1}(X \times \tau) \big),
\end{equation}
where $\tau < \Delta^k$, $|\tau| = k-1$.

\begin{figure}[!htp]
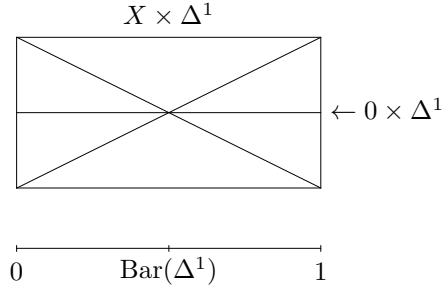

\[
\xy 0;<2cm,0cm>: (-1,0) ; (1,0) **@{-}, (-1,1)
; (1,1) **@{-}, (-1,0) ; (-1,1) **@{-}, (1,0) ; (1,1) **@{-}, (-1,0.5) ;
(1,0.5) **@{-}, (-1,0) ; (1,1) **@{-}, (-1,1) ; (1,0)
**@{-},  (-1,-0.4) ; (1,-0.4) **@{-}, (-1,-0.42) ; (-1,-0.38) **@{-},
(0,-0.42) ; (0,-0.38) **@{-}, (1,-0.42) ; (1,-0.38) **@{-},
(1.35,0.52)*={\quad \leftarrow 0 \times \Delta^1};
(0,-0.55)*={\Bar(\Delta^1)}; (-1,-0.55)*={0}; (1,-0.55)*={1};
(0,1.15)*={X \times \Delta^1}
\endxy
\]
\caption{Picture of $C_1 (X \times \Delta^1)$ for $X = S^0$.}
\label{fig:C-1M}
\end{figure}

Figure \ref{fig:C-1M} depicts the special case when $X = S^0$ and
$k=1$. (Here mutually distinct cone points should be used for
different cones in sight.)

The spaces $CX \times \Delta^k$ and $C_k (X \times \Delta^k)$ can be
identified by choosing suitable cone points in the above
description. Namely, if $\tau \leq \Delta^k$, then choose the cone
point of $C_{|\tau|} (X \times \tau)$ to be the
barycenter $\Bar(\tau)$. Explicitly, we obtain a homeomorphism $C_{|\tau|} (X
\times \tau) \rightarrow CX \times \tau$ as follows.
Suppose that for all $\sigma < \tau$ a homeomorphism $C_{|\sigma|} (
X \times \sigma) \rightarrow CX \times \sigma$ was constructed. Then
define a homeomorphism from $C(X \times \tau \cup CX \times \partial
\tau)$ to $CX \times \tau$ by
\begin{equation*}
\begin{cases}
r(x,t)+(1-r)\Bar(\tau) \mapsto (rx,rt+(1-r)\Bar(\tau)), \textrm{ if
} (x,t) \in X \times \tau,
\\ r(sx,t)+(1-r)\Bar(\tau) \mapsto
(rsx,rt+(1-r)\Bar(\tau)), \textrm{ if } (sx,t) \in CX \times
\partial \tau.
\end{cases}
\end{equation*}
By induction over $\tau \leq \Delta^k$ we obtain the desired
identification. Again, Figure \ref{fig:C-1M} in the special case when
$X = S^0$ and $k=1$ might be helpful.

With this identification, starting with a block map $h \colon X
\times \Delta^k \rightarrow Y \times \Delta^k$ we construct a block
map $C_k (h) \colon CX \times \Delta^k \rightarrow CY \times
\Delta^k$ by iterating the construction of the cone map as described
in the case $k=0$. The resulting map clearly has the property
$C_k(h) (r \cdot X \times \Delta^k) \subseteq r \cdot Y \times
\Delta^k$; and it is our definition of the map $C^+ (h)$ on $CX
\times \Delta^k$. It is easy to verify that the maps $T'$, $f'$
defined via this procedure define a structure on $V \times \Delta^k$
as required.
\end{con}

\begin{proof}[Proof of Proposition \ref{spheres-vs-vector-spaces}]
The composition $\res \circ \cone$ is the identity. It remains to
show that the composition $\cone \circ \res$ is homotopic to the
identity. We start with a structure $(T,f)$ on $V \times \Delta^k$,
consider its restriction, the structure $(T|,f|)$ on $S(V) \times
\Delta^k$ and the image $(T',f') = \cone (T,f)$, a structure on $V
\times \Delta^k$. We need to find compatible homotopies
between $T$ and $T'$ and $f$ and $f'$ satisfying certain
restrictions .

Similarly as in the previous construction we give a general
procedure as follows. Suppose $h \colon C^+ X \times \Delta^k
\rightarrow C^+Y \times \Delta^k$ is a block map such that $h (r \cdot
X \times \Delta^k) \subseteq (r \cdot Y \times \Delta^k)$. Consider
the restriction block map $h| \colon X \times \Delta^k \rightarrow Y
\times \Delta^k$ and the block map $C^+_k (h|) \colon C^+ X \times
\Delta^k \rightarrow C^+ Y \times \Delta^k$ obtained via the
iterated cone construction. We give a construction of a homotopy,
say $h_t$, between $h$ and $C^+_k (h|)$ with the property that each
$h_t$ is a block map $C^+X \times \Delta^k \rightarrow C^+ Y \times
\Delta^k$ such that $h_t (r \cdot X \times \Delta^k) \subseteq (r
\cdot Y \times \Delta^k)$. Applying this homotopy to the maps $T$,
$f$ in place of $h$ produces the desired homotopies between $T$,
$f$ and $T'$, $f'$.

\medskip\noindent \textbf{The case $k=0$.} The required homotopy is
obtained using the Alexander trick. We can think of the map $h$ as
of a perturbation of $C^+ (h|)$ and refer to the part of $C^+ X$ on
which the two maps $h$ and $C^+ (h|)$ do not agree as the domain
of the perturbation. Further we consider $C^+ X = CX \cup C^{1+} X$
and we give a different description of the homotopy on the two pieces. In
words, the trick consists of shrinking the domain of the
perturbation in $CX$ and of pushing the domain of the perturbation
towards $\infty$ in $C^{1+} X$.

Formally, the homotopy $h_t$ on $CX$ is parametrized by $t \in
[0,1]$ and given by
\begin{equation}
h_t (rx) = \begin{cases} rh(1x) & \textrm{ for } r\geq t \\
rh(x,r/t) & \textrm{ for } r < t.
\end{cases}
\end{equation}
The homotopy $h_t$ on $C^{1+} X$ is parametrized by $[1,\infty]$ and
given by
\begin{equation}
h_t (rx) = \begin{cases} rh(1x) & \textrm{ for } r \leq t
\\ rh((1+r-t)x) & \textrm{ for } r \geq t,
\end{cases}
\end{equation}
if $t \in [1,\infty)$ and by $h_\infty (rx) = rh(1x)$. After
reparametrizing, using a suitable homeomorphism $[1,\infty] \cong
[0,1]$ which maps $1$ to $1$ and $\infty$ to $0$, we obtain that
$h_0 = C^+ (h|_X)$ and $h_1 = h$. Moreover, clearly $h_t (r \cdot X)
\subseteq r \cdot Y$ for all $t \in [0,1]$.

\medskip\noindent\textbf{The cases $k>0$.} Again, the homotopy is described on the two
pieces of the decomposition $C^+ X \times \Delta^k = CX \times
\Delta^k \cup C^{1+}X \times \Delta^k$ differently. On $C^{1+}X
\times \Delta^k$ it is given by pushing the domain of the
perturbation towards $\infty$ just as in the case $k=0$; the formula
is left for the reader.

To describe the homotopy on $CX \times \Delta^k$, recall that in the
definition of $C_k(h|)$ we iterated the cone construction. Further,
in the case $k=0$ the homotopy between $h$ and $C(h|)$ was given by
the Alexander trick. In analogy, in the case $k>0$ the homotopy
between $h$ and $C_k(h|)$ is given by iterating the Alexander trick.

To start, recall the identification of $C_k (X \times \Delta^k)$
with $CX \times \Delta^k$ and apply the Alexander trick homotopy to
$h$ considered as a block map from $C_k (X \times \Delta^k)$ to $C_k
(Y \times \Delta^k)$. The resulting map, say $h_1$, coincides with
$C_k (h|)$ on $C(X \times \Delta^k)$, but is possibly different on
the complement, the union of $C(C_{k-1}(X \times \tau))$ where $\tau
< \Delta^k$, $|\tau| = k-1$ (see Figure \ref{fig:C-1M}). Note
however, that if we consider $C(C_{k-1}(X \times \tau))$ for each
such $\tau$ as the quotient of $C_{k-1} (X \times \tau) \times
[0,1]$, then the restriction of $h_1$ is the quotient map of the
product map of the restriction of $h$ to $C_{k-1} (X \times \tau)$
with the identity on $[0,1]$ (again see Figure \ref{fig:C-1M}).
Hence, we can perform a homotopy given as the product of the
Alexander trick homotopy on $C_{k-1} (X \times \tau)$ and the
identity on $[0,1]$. The resulting map now agrees with $C_k(h|)$ on
a larger portion of $CX \times \Delta^k$. We proceed inductively
until we take care of all the faces of $\Delta^k$. Then the
resulting map agrees with $C_k(h|)$ on the whole $CX \times \Delta^k$.
Again, it is clear that the homotopy $h_t$ is such that for all $t$ we have $h_t (r \cdot X \times \Delta^k) \subset r \cdot Y \times \Delta^k$ for all $r \geq 0$.
\end{proof}

Before we discuss continuity of the functor $F$ recall that, since
the maps $T$, $f$ can vary continuously, the set of $k$-simplices of
the $\Delta$-groupoid $\Stwo(V)$ (similarly $\Stwo(S(V))$) is
endowed with a topology, namely the subspace topology of the product
topology of the compact-open topology on the space of self-maps on
$V \times \Delta^k$. So we have a topological $\Delta$-groupoid, and the
associated $\Delta$-set $\Stwo(V)$ becomes a $\Delta$-space.

\begin{lem}
The geometric realization of $\Stwo(V)$ as a $\Delta$-space is
homotopy equivalent to the geometric realization of the same object
as a $\Delta$-set.
\end{lem}

\begin{proof}
To see this we resolve the $\Delta$-space version of $\Stwo(V)$ by
taking the singular simplicial set of the space of $k$-simplices for
each $k \in \N$. We obtain a combined simplicial-$\Delta$-set, which
has its geometric realization homotopy equivalent to the geometric
realization of the $\Delta$-space. We can `change the order' of
realization and realize first spaces $\maps (\Delta^k, \Stwo(V))$
for each $k \in \N$ as $\Delta$-sets and then realize the resulting
simplicial set. The lemma will be proved if we show that all
degeneracy and face maps between $\maps (\Delta^k, \Stwo(V))$ and
$\maps (\Delta^{k+1}, \Stwo(V))$ induce homotopy equivalences. This
can be verified using the Kan property of $\Stwo(V)$.
\end{proof}

Now we show that the functor $F$ is a continuous functor in the
sense of Definition \ref{continuous}. First, note that one can
define a continuous functor from $\sJ$ to topological
$\Delta$-groupoids or to $\Delta$-spaces to be a functor for which
for all $k$ the evaluation maps
\[
\mathrm{ev}_k \colon \mor(V,W) \times E(V)_k \rightarrow E(W)_k
\]
are continuous. Given such a continuous functor from $E \colon \sJ
\rightarrow \mathrm{top-}\Delta\mathrm{-groupoids}_\ast$, composing
with the functor $E \mapsto E^\sharp$ and taking the geometric
realization of $E^\sharp$ we obtain a continuous functor to
$\Spaces_\ast$.

\begin{thm}\label{F-cont}
The functor $F \colon \sJ \rightarrow \Spaces_\ast$ is a continuous
functor.
\end{thm}

\begin{proof}
In view of the previous remarks it is enough to show that the
functor $F \colon V \mapsto \Stwo (V)$ is continuous. That means we
need to show that the evaluation map
\[
\mathrm{ev}_k \colon \mor (V,W) \times F(V)_k \rightarrow F(W)_k
\]
is continuous for each $k$.

If $\dim  (V) = \dim (W)$, then $\mathrm{ev}_k$ is given by
$(\xi,(T,f)) \mapsto (\xi \circ T \circ \xi^{-1}, \xi \! \circ \! f
\! \circ \! \xi^{-1})$ which clearly gives continuity.

For the case when $\dim (V) < \dim (W)$ we work locally and reduce
to the above case as follows. Choose $\xi \colon V \rightarrow W$.
We prove that the evaluation map is continuous on some neighborhood
$V_\xi$ of $\xi$. There is a map $\mor (W,W) \rightarrow \mor (V,W)$
given by $\zeta \mapsto \zeta \! \circ \! \xi$, which is a bundle
map. Therefore there is a neighborhood $V_\xi \subset \mor (V,W)$ of
$\xi$ over which we have a section $s$ of this bundle. Then the
evaluation map $\mathrm{ev}_k$ can be factored as
\[
V_\xi \times F(V)_k \xrightarrow{s \times F(\xi)_k} \mor (W,W)
\times F(W)_k \xrightarrow{\mathrm{ev}_k} F(W)_k.
\]
The first map is continuous because by the assumption both of its
factors are continuous and the second one is the previous case.
\end{proof}

%%%%%%%%%%%%%%%%%%%%%%%%%%%%%%%%%%%%%%%%%%%%%%%%%%%%%%%%%%%%%%
%CORRECTION%
%%%%%%%%%%%%%%%%%%%%%%%%%%%%%%%%%%%%%%%%%%%%%%%%%%%%%%%%%%%%%%

%%%%%%%%%%%%%%%%%%%%%%%%%%%%%%%%%%%%%%%%%%%%%%%%%%%%%%%%%%%%%%%%%%%

\section{Orthogonal Calculus} \label{chap:oc}

Our method of study of the functor $F$ defined in the previous
section is the orthogonal calculus of functors. In this section we
review some of the terms and tools of orthogonal calculus which we
use in this paper. It is intended as a short survey which describes
the context in which we work. For the full account we refer the
reader to the paper of Weiss \cite{Weiss} in which orthogonal
calculus was developed.

We recall some notation. By $\sJ$ is denoted the category of
finite-dimensional real vector spaces with (positive definite) inner
product. Morphisms are linear maps preserving inner product, in
particular they are linear inclusions. The objects of $\sJ$ are
usually denoted by letters $U, V, W$. Morphism sets $\mor(V,W)$ are
Stiefel manifolds and so have a topology. By $\iota$ will be denoted
the inclusion of $V$ into $\VU$ for any objects $V$, $U$ in $\sJ$.
In special cases, by $\iota_1$ will be denoted the inclusion of
$\VR$ into $\VR \oplus \R$ where $\R$ is mapped isomorphically onto
$\R \oplus \{0\}$ in $\R \oplus \R$ and by $\iota_2$ the inclusion
where $\R$ is mapped onto $\{0\} \oplus \R$ in $\R \oplus \R$. For
$U \in \mathrm{Ob}(\sJ)$, $U^c$ denotes the one-point
compactification of $U$. Furthermore $\R_+ = \{ u \in \R \; | \; u
\geq 0 \}$ and $\R_+^c = \R_+ \cup \{ \infty \}$. We also write $k
\cdot U$ for $\R^k \otimes U$.

The orthogonal calculus works for functors $E \colon \sJ \rightarrow
\Spaces_\ast$ satisfying the following axiom.

\begin{defin} \cite{Weiss} \label{continuous}
A functor $E \colon \sJ \rightarrow \Spaces_\ast$ is called {\it
continuous} if for each $V, W \in \mathrm{Ob}(\sJ)$ the evaluation
map
\[ \mathrm{ev} \colon \mor (V,W) \times E(V) \rightarrow E(W) \] is
continuous.
\end{defin}

We have shown in Theorem \ref{F-cont} that the functor $F$ satisfies
this axiom. There are many other functors of this type. For example,
the functor $V \mapsto \mathrm{BO}(V)$ where $\mathrm{BO}(V)$ is the
classifying space of the orthogonal group $\mathrm{O}(V)$ has been
studied in \cite{Arone}. Other examples come up in the study of
automorphisms of manifolds, see papers \cite{Weiss}, \cite{WW}.

We have already presented the basic philosophy of orthogonal
calculus in the introduction. Recall that a continuous functor $E
\colon \sJ \rightarrow \Spaces_\ast$ is studied by means of its
Taylor tower which is a tower of functors
\begin{equation} \label{Taylor}
\cdots \rightarrow T_kE \rightarrow T_{k-1}E \rightarrow \cdots
\rightarrow T_0E,
\end{equation}
where the functor $T_kE$ is the $k$-th polynomial approximation of
$E$. The tower should be thought of as the expansion at infinity.
There are two aspects of the Taylor tower which one needs to
understand. Firstly, one needs to find out whether the Taylor tower
converges to $E$, that means to estimate the connectivity of the
canonical maps $E(V) \rightarrow T_kE(V)$. Secondly, one needs to
determine the differences, $\hofiber (T_kE(V) \rightarrow
T_{k-1}E(V))$, between the stages of the tower.

Both issues can be investigated using the derivatives in orthogonal
calculus. These arise as follows. Given a continuous functor $E
\colon \sJ \rightarrow \Spaces_\ast$ there is defined for each $k
\in \N$ its $k$-th {\it derivative functor} $E^{(k)}$ which is again
a continuous functor from $\sJ$ to $\Spaces_\ast$ with the
additional structure of a coordinate free spectrum of multiplicity
$k$. That means there are maps
\[ \sigma_k \colon (k \cdot U)^c \wedge E^{(k)} (V) \rightarrow
E^{(k)}( \VU).
\]
On objects $E^{(k+1)}$ is defined inductively as
\[ E^{(k+1)} (V) = \hofiber (\ssigma_k \colon E^{(k)}(V) \rightarrow \Omega^k
E^{(k)}(\VR)),
\]
where $\ssigma_k$ is the adjoint of $\sigma_k$. For the definition
on morphisms, take $\xi \! \in \! \mor(V,W)$ and define $\xi' \in
\mor(\VR,W \oplus \R)$ by $\xi'(v,x) = (\xi(v),x)$. The pair
$(\xi,\xi')$ then defines the map $E^{(k+1)}(\xi)$ between the
homotopy fibers $E^{(k+1)}(V) \rightarrow E^{(k+1)}(W)$. In
Construction \ref{map1} we confine ourselves to a description of the
maps $\sigma_k$ for $k=1$ as we will not need the maps $\sigma_k$
for $k > 1$.

\begin{con} \label{map1}
The values $E^{(1)}(V)$ are defined as homotopy fibers of the
stabilization maps $E(\iota) \colon E(V) \rightarrow E(V \oplus \R)$
where $\iota \colon V \rightarrow V \oplus \R$ is the linear
inclusion. A point in $E^{(1)}(V)$ is described as a pair $(x,
\lambda)$, where $x \in E(V)$ and $\lambda \colon \R_+^c \rightarrow
E(\VR)$ satisfies $\lambda(0) = E(\iota)(x)$ and $\lambda(\infty) =
\ast$.

We give a description of the adjoints of the maps $\sigma_1$, that
means the maps
\[
\ssigma_1 \colon E^{(1)}(V) \rightarrow \Omega^U E^{(1)}(V \oplus
U).
\]
Given $(x,\lambda) \in E^{(1)}(V)$ we describe the image
$(x',\lambda') := \ssigma (x,\lambda) \in \Omega^U E^{(1)}(V \oplus
U)$ where
\[
x' \colon U^c \rightarrow E(V \oplus U),
\]
\[
\lambda' \colon (U \oplus \R_+)^c \rightarrow E(V \oplus U \oplus
\R),
\]
such that $x'(\infty) = \ast$, $\lambda'|_{(U \oplus \{0\})^c} =
E(\iota)(x')$ and $\lambda' (\infty) = \ast$.

Each vector $u \in U$, such that $\|u\| = 1$, defines a morphism
$\xi_u \colon \VR \rightarrow V \oplus U$ in $\sJ$ by $\xi_u \colon
(v,1) \mapsto (v,u)$.

Similarly, each vector $w \in U \oplus \R_+$, such that $\|w\| = 1$,
defines a morphism $\zeta_w \colon \VR \rightarrow V \oplus U \oplus
\R$ in $\sJ$ by $\zeta_w \colon (v,1) \mapsto (v,w)$.

All such unit vectors $u \in U$, $w \in U \oplus \R_+$ define rays
$\langle u \rangle$ and $\langle w \rangle$ from $0$ trough $u$ and
$w$, which cover $U$ and $U \oplus \R_+$ respectively. We define
$x', \lambda'$ on the closures $\langle u \rangle^c \subseteq U^c$,
$\langle w \rangle^c \subseteq (U \oplus \R_+)^c$ as the
compositions
\[
x' \colon \langle u \rangle^c \xrightarrow{\cong} \R_+^c
\xrightarrow{\lambda} E(\VR) \xrightarrow{E(\xi_u)} E(V \oplus U),
\]
\[
\lambda' \colon \langle w \rangle^c \xrightarrow{\cong} \R_+^c
\xrightarrow{\lambda} E(\VR) \xrightarrow{E(\zeta_w)} E(V \oplus U
\oplus \R).
\]
Here the fact that the functor $E$ is continuous is used. It
guarantees that the maps $x'$, $\lambda'$ are continuous.
\begin{figure}[!htp]
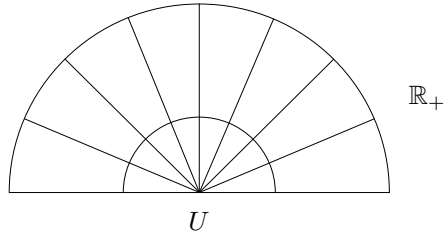

\[
\xy 0;<2.5cm,0cm>: (0,0)*\cir<2.5cm>{u^d}, (0,0)*\cir<1.cm>{u^d},
,a(0) **\dir{-} ,a(23) **\dir{-} ,a(45) **\dir{-} ,a(67) **\dir{-}
,a(90) **\dir{-} ,a(112) **\dir{-} ,a(135) **\dir{-} ,a(157)
**\dir{-} ,a(180) **\dir{-}
;(0,-0.15)*={U};(1.2,0.5)*={\R_+}
\endxy
\]
\caption{Schematic picture of the map $\ssigma_1$}
\label{fig:sigma1}
\end{figure}

Figure \ref{fig:sigma1} describes the map $\ssigma_1$. Lines inside
the half-disk represent the rays $\langle w \rangle$ and the
half-circle represents the unit hemisphere inside $U \oplus \R_+$.
The morphism $\zeta_w \in \mor(\VU,\VU \oplus \R)$ which is used to
define $\lambda'$ on each ray $\langle w \rangle$ is determined by
the point on the unit hemisphere which the ray $\langle w \rangle$
intersects.
\end{con}

Now we turn our attention to the convergence question. We will not
address this issue in full generality, that means we will not
discuss a notion of an analytic functor analogous to such a notion
in \cite{Go}. We will only consider the case when a functor $E$ is
polynomial of some degree $\leq k$. In this case the Taylor tower
`stops' at the $k$-stage, that means the canonical maps $E(V)
\rightarrow T_kE(V)$ are homotopy equivalences and thus the $k$-th
stage $T_kE$ gives us complete homotopy theoretical information
about the functor $E$. Polynomial functors are defined as follows.

\begin{defin}
A continuous functor $E \colon \sJ \rightarrow \Spaces_\ast$ is
called {\it polynomial of degree} $\leq k$ if for every $V \in \sJ$
the following canonical map is a homotopy equivalence
\[ E(V) \longrightarrow \holimsub{0 \neq U \subseteq \R^{k+1}} E(\VU). \]
The homotopy limit is taken over the topological poset of all
non-zero vector subspaces of $\R^{k+1}$ (see \cite[section
5]{Weiss}).
\end{defin}

In general, for a functor $E \colon \sC \rightarrow \Spaces_\ast$
from a small category $\sC$, the homotopy limit of $E$ over $\sC$
can be defined as the total space of a cosimplicial space given by
\begin{equation} \label{cosimpl-space}
\l \mapsto \prod_{\sigma \colon \l \rightarrow \sC} E(\sigma(l)),
\end{equation}
see \cite{BK}. If $\sC$ is a topological category, this affects the
definition of the homotopy limit. In our case, when $\sC$ is $\{ 0
\neq U \subseteq \R^{k+1} \}$ with objects non-zero vector subspaces
of $\R^{k+1}$ and one morphism for every inclusion $U \subseteq V
\subseteq \R^{k+1}$, the sets $\{ \sigma \colon \l \rightarrow \sC
\}$ become disjoint union of Stiefel manifolds. Then the disjoint
union of spaces $E(\sigma(l))$ becomes a bundle over $\{ \sigma
\colon \l \rightarrow \sC \}$. The space of $l$-cosimplices of the
cosimplicial space which defines the homotopy limit then becomes the
space of sections of this bundle, see \cite[sections 4,5]{Weiss} for
more details.

The property above is an extrapolation property. Notice that the
value of the functor $E$ at $V$ is determined by the values of $E$
at vector spaces between $V$ and $\VR^{k+1}$ which are strictly
greater than $V$. For more details see \cite[section 5]{Weiss}.

The derivatives are related to the polynomial functors via the
following proposition \cite[Proposition 5.3.]{Weiss}.

\begin{prop} \label{fibseq}
For any $V \in \sJ$ the value $E^{(k+1)}(V)$ is homotopy equivalent
to the homotopy fiber of the canonical map
\[E(V) \rightarrow \holimsub{0 \neq U \subseteq \R^{k+1}} E(\VU)\]
over the base point of the target.
\end{prop}

Thus if the functor $E$ is polynomial of degree $\leq k$ then the
$(k+1)$-th derivative functor $E^{(k+1)}$ vanishes. If the target of
the map from Proposition \ref{fibseq} is path-connected for all $V$,
then the converse is also true. If the target is not path-connected,
then we have to investigate also the homotopy fibers over other
points than the base point. This leads to a necessary and sufficient
condition for $E$ to be polynomial of degree $\leq k$ which we state
in Theorem \ref{poly1} in the case $k=1$. It is due to Weiss (but
did not appear in \cite{Weiss}) and we are grateful to him for
showing it to us. Before we state it we need to define a
modification of the map $\ssigma_1$ presented in Construction
\ref{map2}.

We also need some new notation because we will use compactifications
of Euclidean spaces other than the one-point compactification. For
$V \in \mathrm{Ob}(\sJ)$ by $\overline{V}$ is denoted the
compactification of $V$ by the sphere $\partial \overline{V} =
\{\overline{w} \; | \; w \in V, \; \|w\|=1 \}$ at infinity. In
special cases, $\Rbar$ is the $0$-sphere compactification of $\R$,
with the points at infinity $\pm \overline{1}$, $\overline{\R \oplus
\R}$ is the $1$-sphere compactification of $\R \oplus \R$ with
points at infinity $\partial \; \overline{\R \oplus \R} = \{
\overline{w} \; | \; w \in \R^2, \; \| w \| = 1 \}$. Here we also
distinguish $\partial_\pm \overline {\R \oplus \R} = \{ \overline{w}
\; | \; w=(w_1,w_2), \, \pm w_2 \geq 0 \}$. By $\Rtwoplusbar$ is
denoted the $1$-disk compactification of $\R \oplus \R_+$, points at
infinity are $\{ \wbar \; | \; w \in \R \oplus \R_+, \; \|w\|=1 \}$.

\begin{con} \label{map2}
Denote $E^{(1)}(V,y) = \hofibery (E(V) \rightarrow E(\VR))$, where
$\hofibery$ means the homotopy fiber taken over some point $y \in
E(\VR)$.

Let $X_y$ denote the space of pairs of maps $(x',\lambda')$ where
\[
x' \colon \Rbar \rightarrow E(\VR),
\]
\[
\lambda' \colon \Rtwoplusbar \rightarrow E(\VR^2),
\]
such that $x'(\pm \overline{1}) = \pm y$, where $-y$ denotes the
image of the point $y$ under the self-map of $E(\VR)$ induced by the
reflection of $\R$. Furthermore we require that $\lambda'|_{\Rbar
\oplus \{0\}} = E(\iota_1)(x')$ and $\lambda' (\wbar) =
E(\zeta_w)(y)$ where $\zeta_w$ is as in Construction \ref{map1}.

Notice that in the case $y = \ast$ the space $X_y$ is homotopy
equivalent to the loop space $\Omega E^{(1)}(\VR)$. If $y \neq \ast$
the space $X_y$ can be thought of as a certain space of paths in
$E^{(1)}(\VR,y')$ where $y' = E(\zeta_{w'})(y)$ where $w' =
(0,0,1)$. This is depicted in Figure \ref{fig:sigma2}.

\begin{figure}[!htp]
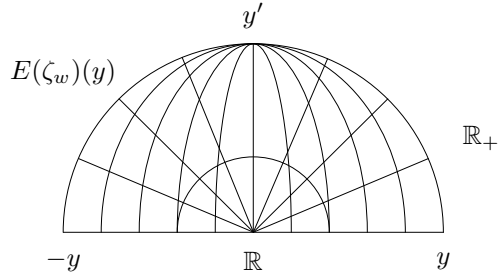

\[
\xy 0;<2.5cm,0cm>: (0,0)*\cir<1.cm>{u^d}, (0,1);**\dir{-}, (-1,0);
{\ellipse ,_v(1,0){}}; (-0.8,0); {\ellipse ,_v(1,0){}}; (-0.6,0);
{\ellipse ,_v(1,0){}}; (-0.4,0); {\ellipse ,_v(1,0){}}; (-0.2,0);
{\ellipse ,_v(1,0){}}; (0.2,0); {\ellipse ,v(1,0){}}; (0.4,0);
{\ellipse ,v(1,0){}}; (0.6,0); {\ellipse ,v(1,0){}}; (0.8,0);
{\ellipse ,v(1,0){}}; (1,0); {\ellipse ,v(1,0){}}; (-1,0) ; (1,0)
**@{-}; (0,0); a(23) **\dir{-}, a(45) **\dir{-}, a(67) **\dir{-},
a(112) **\dir{-}, a(135) **\dir{-}, a(157) **\dir{-},
(0,-0.15)*={\R}; (1.2,0.5)*={\R_+}; (-1,-0.15)*={-y};
(1,-0.15)*={y}; (0,1.15)*={y'}; (-1,0.85)*={E(\zeta_w)(y)}
\endxy
\]
\caption{Schematic picture of $\lambda'$ and of the map
$\ssigma_1$.} \label{fig:sigma2}
\end{figure}

The map $\ssigma_1 \colon E^{(1)}(V,y) \rightarrow X_y$ is defined
similarly as the map $\ssigma_1$ in Construction \ref{map1}. For
$(x,\lambda) \in E^{(1)}(V,y)$ we define $(x',\lambda') = \ssigma_1
(x,\lambda) \in X_y$ as follows. Given unit vectors $u \in \R$, $w
\in \R \oplus \R_+$ we now have rays $\langle 0,\ubar \rangle$,
$\langle 0,\wbar \rangle$ from $0$ to $\ubar$ or $\wbar$. These rays
cover $\Rbar$, $\Rtwoplusbar$ respectively. Define
\[
x' \colon \langle 0,\ubar \rangle \xrightarrow{\cong} \R_+^c
\xrightarrow{\lambda} E(\VR) \xrightarrow{E(\xi_u)} E(\VR),
\]
\[
\lambda' \colon \langle 0,\wbar \rangle \xrightarrow{\cong} \R_+^c
\xrightarrow{\lambda} E(\VR) \xrightarrow{E(\zeta_w)} E(\VR^2).
\]
If $y = \ast$ then the map $\ssigma_1$ obviously coincides with the
map $\ssigma_1$ from Construction \ref{map1}, see Figure
\ref{fig:sigma2} again.
\end{con}

Now we are ready to state the promised necessary and sufficient
condition.

\begin{thm}\label{poly1}
A continuous functor $E \colon \sJ \rightarrow \Spaces_\ast$ is
polynomial of degree $\leq 1$ if and only if for all $V \in \sJ$ and
for all $y \in E(\VR)$ the map $\ssigma_1 \colon E^{(1)}(V,y)
\rightarrow X_y$ defined in \ref{map2} is a homotopy equivalence.
\end{thm}

\begin{proof}

We sketch the proof. We need to show that the map
\begin{equation} \label{poly}
E(V) \rightarrow \holimsub{0 \neq U \subset \R^2} E(\VU)
\end{equation}
is a homotopy equivalence if and only if for all choices of $V$ and
$y$ the map $\ssigma_1$ is.

Consider the map $E(V) \rightarrow E(\VR)$ induced by $\iota \colon
V \rightarrow \VR$ and the map $\mathrm{holim}_{0 \neq U \subset
\R^2} E(\VU) \rightarrow E(\VR)$ obtained by the restriction via
$\iota_1 \colon \VR \rightarrow \VR^2$. For $y \in E(\VR)$ consider
the induced map, say $\ssigma_y$, between the homotopy fibers over
$y$ of both maps. The second map is a fibration so the target
homotopy fiber can be replaced by the fiber over $y$. In Lemma
\ref{X_y} below we show that the space $X_y$ is homotopy equivalent
to this fiber. The map $\ssigma_1$ is such that it is homotopic to
the composition of the map $\ssigma_y$, a suitable homotopy
equivalence from the  homotopy fiber to the fiber over $y$ of the
map $\mathrm{holim}_{0 \neq U \subset \R^2} E(\VU) \rightarrow
E(\VR)$ and the identification of this fiber with the space $X_y$
given in Lemma \ref{X_y} (this should be clear after reading the
proof of Lemma \ref{X_y} below). Now the map (\ref{poly}) is a
homotopy equivalence if and only if for all choices of $y$ the map
$\ssigma_y$ is.
\end{proof}

\begin{lem} \label{X_y}
The space $X_y$ is homotopy equivalent to the fiber of the map
\[
\holimsub{0 \neq U \subset \R^2} E(\VU) \rightarrow E(\VR)
\]
over the point $y$.
\end{lem}

\begin{proof}
As a topological space the category $\{0 \neq U \subseteq \R^2 \}$
is $\R P(\R^2) \coprod \{\R^2\}$, with one morphism (an inclusion,
denoted by $\iota_U $) for every point of $U \in \R P(\R^2)$ to
$\R^2$. We will denote by the same $\iota_U$ also the inclusion of
$\VU$ into $\VR^2$.

According to the definition the homotopy limit of the functor $E(V
\oplus -)$ over this category is given by a pair of maps
\begin{eqnarray*}
\lambda_0 & \colon & \R P(\R^2) \rightarrow E_0 \\ \lambda_1 &
\colon & \R P(\R^2) \times \Delta^1 \rightarrow E(\VR^2)
\end{eqnarray*}
where $E_0$ is a bundle over the circle $\R P(\R^2)$ with the fiber
over $U$ the space $E(\VU)$ and such that $\lambda_0$ is a section
of this bundle, $\lambda_1 (U,0) = E(\iota_U) (\lambda_0 (U))$ and
$\lambda_1 (U,1) = \lambda_1 (U',1)$ for every $U$, $U'$ in $\R
P(\R^2)$. In a more condensed form we can say that the source of the
map $\lambda_1$ is defined just on the cone over the circle $\R
P(\R^2)$. The projection to $E(\VR)$ is given by $\lambda_0 (\R
\oplus \{0\})$, therefore taking the fiber over $y \in E(\VR)$ means
to require in addition that $\lambda_0 (\R \oplus \{0\}) = y$.

Notice that we can think of $E_0$ as a subbundle of the trivial
bundle over $\R P(\R^2)$ with the fiber $E(\VR^2)$. The
compatibility conditions between $\lambda_0$ and $\lambda_1$ express
the fact that the two sections of the two bundles commute with the
embedding of the subbundle. The identification of the fiber over $y$
with the space $X_y$ is obtained via a sequence of
reparametrizations of the maps $\lambda_0$, $\lambda_1$ which
correspond to certain isotopies of the embedding of $E_0$ into
$E(\VR^2) \times \R P(\R^2)$ as follows. (The reader might find it
useful to draw some pictures of the situation similar to Figure
\ref{fig:sigma2}.)

Firstly, recall the $1$-sphere compactification $\overline{\R \oplus
\R}$ of $\R \oplus \R$ with points at infinity $\partial \;
\overline{\R \oplus \R} = \{ \overline{w} \; | \; w \in \R^2, \; \|
w \| = 1 \}$ and redefine the fiber over $y$ as the space of maps
\begin{eqnarray*}
\lambda_0 & \colon & \partial \; \overline{\R \oplus \R} \rightarrow
E_0
\\ \lambda_1 & \colon & \overline{\R \oplus \R} \rightarrow
E(\VR^2)
\end{eqnarray*}
such that for $w = (w_1,w_2) \in \R \oplus \R$ with $\|w\| = 1$ we
have $\lambda_0 (\overline{w}) \in E(V \oplus \langle w \rangle)$ if
$w_2 \leq 0$ and $\lambda_0 (\overline{w}) = y \in E(V \oplus \R
\oplus \{0\})$ if $w_2 \geq 0$, and $\lambda_1 (\overline{w}) = E
(\iota_{\langle w \rangle}) (\lambda_0 (\overline{w}))$ if $w_2 \leq
0$ and $\lambda_1 (\overline{w}) = E(\iota_1) (y)$ if $w_2 \geq 0$.
Here $\langle w \rangle$ denotes the vector subspace of $\R^2$
generated by $w$.

By `rotating' the fibers of $E_0$ inside $E(\VR^2)$ redefine our
space again as the space of maps $(\lambda_0,\lambda_1)$ as above,
but satisfying now $\lambda_0 (\overline{w}) \in E(V \oplus \R)$ if
$w_2 \leq 0$ and $\lambda_0 (\overline{w}) = E(\xi_w) (y) \in E(V
\oplus \langle w \rangle)$ if $w_2 \geq 0$, and $\lambda_1
(\overline{w}) = E(\iota_1) (\lambda_0 (\wbar))$ if $w_2 \leq 0$ and
$\lambda_1 (\overline{w}) = E (\iota_{\langle w \rangle}) (\lambda_0
(\overline{w})) = E(\zeta_w) (y)$ if $w_2 \geq 0$, where $\xi_w
\colon \VR \rightarrow V \oplus \langle w \rangle$ and $\zeta_w
\colon \VR \rightarrow \VR^2$ are as in Construction \ref{map1}.

In the next step consider $\partial_- \overline {\R \oplus \R}$ and
redefine the fiber over $y$ as the space of maps
\begin{eqnarray*}
\lambda_0 & \colon & \partial_- \overline{\R \oplus \R} \rightarrow
E (\VR)
\\ \lambda_1 & \colon & \overline{\R \oplus \R} \rightarrow
E(\VR^2)
\end{eqnarray*}
such that for $w = (w_1,w_2) \in \R \oplus \R$ with $\|w\| = 1$ we
have $\lambda_0 (\overline{(\pm 1,0)}) = \pm y$, and $\lambda_1
(\overline{w}) = E (\iota_1) (\lambda_0 (\overline{w}))$ if $w_2
\leq 0$ and $\lambda_1 (\overline{w}) = E(\zeta_w) (y)$.

Finally choose a homeomorphism $\overline{\R \oplus \R} \rightarrow
\Rtwoplusbar$ such that its restriction to $\partial_- \overline {\R
\oplus \R}$ is the projection onto $\Rbar \oplus \{0\}$ and the
restriction to $\partial_+ \overline {\R \oplus \R}$ is the
identity. The resulting space of maps is the desired $X_y$.
\end{proof}

\begin{rem}
\label{pattern} Note that if we want to use Theorem \ref{poly1} to
prove that a certain functor $E$ is polynomial of degree $\leq 1$
there are two possibilities for each $y \in E(\VR)$; the spaces
$E^{(1)}(V,y)$, $X_y$ might be either empty or not. So the proof
splits into two parts. Namely,
\begin{enumerate}
\item show that $E^{(1)}(V,y)$ is empty if and only if $X_y$ is
empty, and
\item if $E^{(1)}(V,y)$ is not empty, show that $\ssigma_1 \colon
E^{(1)}(V,y) \rightarrow X_y$ is a homotopy equivalence.
\end{enumerate}
\end{rem}

This finishes the discussion of polynomial functors. Now we turn to
the Taylor tower. We have already mentioned that for a continuous
functor $E$ the functors $T_kE$ are its polynomial approximations.
This is the content of Theorem 6.3 of \cite{Weiss} which in a
slightly weaker form reads as:

\begin{thm} \label{approx}
Let $E \colon \sJ \rightarrow \Spaces_\ast$ be a continuous functor.
Then for $k \geq 0$ there exists a functor $T_kE \colon \sJ
\rightarrow \Spaces_\ast$ and a natural transformation $\eta_k
\colon E \rightarrow T_k E$ such that:
\begin{enumerate}
\item $T_k E$ is polynomial of degree $\leq k$, \item if $E$
already is polynomial of degree $\leq k$ then the map \[ \eta_k
\colon E(V) \rightarrow T_k E(V)\] is a homotopy equivalence for all
$V$.
\end{enumerate}
\end{thm}

The formula for $T_kE$ is given in the proof of Theorem 6.3 in
\cite{Weiss}. For example, if $k=0$ we have for all $V \in
\mathrm{Ob}(\sJ)$
\[
T_0E(V) = E (\VR^\infty) = \hocolimsub{m \in \N} E(\VR^m).
\]

Now we address the issue how to determine the layers of the Taylor
tower. They are determined by the $k$-th {\it derivative spectrum}
$\TTheta E^{(k)}$ which is the `ordinary' spectrum associated to the
coordinate free spectrum given by the $k$-th derivative functor
$E^{(k)}$. Namely $(\TTheta E^{k})_{km} = E^{(k)}(\R^m)$ and the
structure maps are specializations:
\[
\sigma_k \colon (\R^k)^c \wedge (\TTheta E^{(k)})_{km} \rightarrow
(\TTheta E^{(k)})_{k(m+1)}.
\]
As we have mentioned before the tower is the expansion at infinity
and therefore the spectrum $\TTheta E^{(k)}$ should be thought as
the $k$-th derivative of $E$ at {\it infinity}.

We have also already indicated that there is an action of the
orthogonal group $O(k)$ on the $k$-th derivative spectrum associated
to $E^{(k)}$. The relevant statement in \cite{Weiss} is Proposition
3.1 and simplified it reads as:

\begin{prop} \label{action}
There exists a unique family $\{ \alpha_V \}$ of left actions \[
\alpha_V \colon O(k) \times E^{(k)}(V) \rightarrow E^{(k)}(V) \]
which makes the maps \[ \sigma_k \colon (k \cdot U)^c \wedge
E^{(k)}(V) \rightarrow E^{(k)}(\VU) \] into $O(k)$-maps. Here $O(k)$
acts diagonally on the domain of the map $\sigma_k$, it acts on $(k
\cdot U) = \R^k \otimes U$ because it acts on $\R^k$.
\end{prop}

As we have said before, the polynomial approximations $T_kE$ fit
together to form the {\it Taylor tower} of $E$ which is the diagram
(\ref{Taylor}). The natural transformation $T_kE \rightarrow
T_{k-1}E$ is essentially given as the $(k-1)$-st polynomial
approximation of $T_k E$. The usefulness of the Taylor tower follows
from the fact that the differences between its stages, $\hofiber
(T_kE(V) \rightarrow T_{k-1}E(V))$, which are called the layers of
the tower, can be described using the derivatives. The formula is
given in Theorem 9.1 of \cite{Weiss}. It states:

\begin{thm} \label{tower}
For a continuous functor $E \colon \sJ \rightarrow \Spaces_\ast$,
any $k > 0$ and any $V \in \sJ$ there is the following homotopy
fibration sequence:
\[ \Omega^\infty[((k \cdot V)^c \wedge \TTheta E^{(k)})_{hO(k)}] \rightarrow T_k
E(V) \xrightarrow{r_k} T_{k-1} E(V) \] where the subscript $hO(k)$
denotes the homotopy orbit spectrum. The group $O(k)$ acts
diagonally on the smash product, it acts on $\TTheta E^{(k)}$ by
Proposition \ref{action} and it acts on $k \cdot V = \R^k \otimes V$
because it acts on $\R^k$.
\end{thm}

Once we have established the convergence question and we know the
layers of the Taylor tower we can use the Taylor tower to obtain
some information about the functor $E$. One way how to do this is
suggested in the introduction for our functor $F$. The input for the
Taylor tower of $F$ is given in Theorems \ref{thm1}. and \ref{thm2}.
For examples how the Taylor tower can be used to give information
about other functors see \cite[section 10]{Weiss}.

\begin{rem} \label{category}
The functor $F$ defined in section \ref{chap:suspension} is a
continuous functor $\sJ \rightarrow \Spaces_\ast$ which is obtained
from the functor to the category $\Delta\!-\!\mathrm{sets}_\ast$ (or
$\Delta\!-\!\mathrm{spaces}_\ast$) by the geometric realization.
Although we do not generalize the machinery of orthogonal calculus
to this setting we often use constructions in the category
$\Delta\!-\!\mathrm{sets}_\ast$ (or $\Delta\!-\!
\mathrm{spaces}_\ast$), for example the construction of the homotopy
fiber of a map or other constructions in section \ref{chap:ThetaF1}
in Construction \ref{map3}. After geometric realization these
constructions always commute (at least up to homotopy) with the
corresponding constructions in the category $\Spaces_\ast$ which
justifies their use.
\end{rem}

%%%%%%%%%%%%%%%%%%%%%%%%%%%%%%%%%%%%%%%%%%%%%%%%%%%%%%%

\section{The First Derivative Functor $F^{(1)}$}
\label{chap:F1}

Now we start to apply the machinery of orthogonal calculus to study
the functor $F$ defined in section \ref{chap:suspension}. The aim of
this section is to identify the homotopy type of the spaces
$F^{(1)}(V,y)$. This is done in Proposition \ref{crit-F} where we
give a criterion for the space $F^{(1)}(V,y)$ to be empty and we
show that when it is not empty it is homotopy equivalent to a
certain space $\LN_n(\phi_V)$. The homotopy groups of the space
$\LN_n(\phi_V)$ are the groups $LN_{n+\ast}(\phi_V)$ which are
special cases of the $LS$-groups which are obstruction groups in
codimension $1$ surgery theory. The groups $LN_{n+\ast}(\phi_V)$ are
well known to be isomorphic to the $L$-groups of the trivial group
from `ordinary' surgery theory. Thus we also obtain a calculation of
the homotopy groups of the spaces $F^{(1)}(V,y)$.

Relating the comparison between the block structure spaces of real
projective spaces of different dimensions with codimension $1$
surgery theory is not our original idea. At the level of structure
sets it was used by Browder-Livesay, L\'opez de Medrano and Wall in
order to calculate the structure sets of the real projective spaces
(see \cite{BL}, \cite{LdM}, \cite[chapter 14C]{Wall}). The
statements at the level of block structure spaces which we review
here follow from general theory of \cite[chapter 17A]{Wall} and from
the treatment of codimension $q$ surgery theory in \cite[chapter
7]{Ran1}.

%%%%%%%%%%%%%%%%%%%%%%%%%%%%%%%%%%%%%
%Correction%
%%%%%%%%%%%%%%%%%%%%%%%%%%%%%%%%%%%%%

In this section we use the
definition of the block structure space $\S(X)$ from Construction
\ref{blockdef} and we work on the level of $\Delta$-sets. This is
legal by Remarks \ref{models}, \ref{category}. We also adopt the
convention that an embedding of manifolds will always mean locally
flat embedding. This condition is fulfilled in all cases that we
treat here.

%%%%%%%%%%%%%%%%%%%%%%%%%%%%%%%%%%%%%
%Correction%
%%%%%%%%%%%%%%%%%%%%%%%%%%%%%%%%%%%%%

\subsection{Codimension $q$ surgery theory.} We include a brief
review of the codimension $q$ surgery theory in which we fix the
terminology, the notation and we also introduce some constructions
that will be used in section \ref{chap:ThetaF1}. The reference for
the general theory is \cite[chapter 11]{Wall}. The treatment there
is slightly more general, see remarks by the editor at the beginning
of the chapter. Another useful reference is \cite[chapter 7]{Ran1}.
We do assume however a basic knowledge of `ordinary' surgery theory,
in particular, the algebraic definition of the $L$-groups and the
construction of the surgery obstruction for a degree one normal map
(see \cite[chapters 5,6]{Wall}).

First we need the following construction:

\begin{con} \label{split-decomposition}
Let $M \hookrightarrow X$ be an embedding of manifolds. It induces
the following decomposition of the ambient manifold $X$. The
submanifold $M$ has a normal block bundle $\nu$, with total space
$E(\nu)$ and with the total space of the associated sphere bundle
denoted by $S(\nu)$. The closure of the complement of $E(\nu)$ in
$X$ is denoted by $C$. Thus we have
\[
X = E(\nu) \cup_{S(\nu)} C.
\]
\end{con}

\begin{defin} \label{def_split}
Let $M \hookrightarrow X$ be an embedding of manifolds with the
normal block bundle $\nu$. A simple homotopy equivalence $f \colon Y
\rightarrow X$ from a manifold $Y$ of the same dimension as $X$ is
said to be {\it split along} $M$ w.r.t $\nu$, if $f$ is transverse
to $M$ w.r.t. $\nu$ and the following restrictions are simple
homotopy equivalences:
\begin{itemize}
\item $f|_{M'} \colon M'=f^{-1}(M) \rightarrow M$, \item $f|_{C'}
\colon C'=f^{-1}(C) \rightarrow C$.
\end{itemize}
\end{defin}

\begin{defin} \label{splitprob}
Let $M \hookrightarrow X$ be an embedding of manifolds as above. We
will refer to a simple homotopy equivalence $f \colon Y \rightarrow
X$ as a {\it splitting problem} ({\it along} $M$).
\end{defin}

Given a splitting problem $f$ along $M$ the task is to change $f$ by
a homotopy to a simple homotopy equivalence which is split along
$M$. If $q = (\dim (X) - \dim (M)) \geq 3$ the obstruction theory
reduces to the ordinary surgery theory with the $L$-groups as
obstruction groups (\cite[Theorem 11.3]{Wall}). In cases $q=1,2$ the
obstruction theory is different. Firstly, we can assume that $\nu$
is a vector bundle.  Then given a splitting problem $f$ along $M$
there is an obstruction $\theta(f)$ in a group $LS_n (\Phi)$. The
obstruction vanishes if and only if the problem can be solved. Here
$n = \dim (M)$ and $\Phi$ is the following diagram of fundamental
groupoids:
\[
\xymatrix{
  \pi(S(\nu)) \ar[d] \ar[r] & \pi(C) \ar[d] \\
  \pi(M)  \ar[r] & \pi(X).
}
\]
The groups $LS_n (\Phi)$ depend only on $\Phi$, dimension $n$,
codimension $q$ and the orientation characters of $M$ and $X$. See
\cite[chapter 11]{Wall} for the details of the construction of the
$LS$-groups.

In the special case when $\pi(M) \rightarrow \pi(X)$ is an
isomorphism also the other horizontal morphism in $\Phi$ is an
isomorphism and thus the diagram depends only on the vertical
morphism $\phi \colon \pi(C) \rightarrow \pi(X)$. In this case the
group $LS_n(\Phi)$ is denoted $LN_n(\phi)$. This will be the most
important case for us.

Suppose now we are given an embedding of manifolds $M^n
\hookrightarrow X^{n+q}$ which induces the diagram of fundamental
groupoids $\Phi$. Following general philosophy of \cite[chapter
17A]{Wall} to `spacify' all obstruction groups in surgery theory we
can define a space $\LS_n (\Phi)$ such that $\pi_k (\LS_n(\Phi))
\cong LS_{n+k}(\Phi)$ (see also \cite{Ham}). In an outline, a
$k$-simplex of the space $\LS_n (\Phi)$ is a splitting problem $Z
\rightarrow Y$ of manifold $(k+2)$-ads along an embedding $N
\hookrightarrow Y$ of manifold $(k+2)$-ads, which is solved on one
part of the boundary and comes with a reference map $Y \rightarrow
X$ which respects the decompositions of $Y$ and $X$ given by $N
\hookrightarrow Y$ and $M \hookrightarrow X$ as described in
Construction \ref{split-decomposition}.

We will not give a precise definition of this space in general,
instead we look at the special case of an embedding $M^n
\hookrightarrow E(\xi)$ of a manifold as the zero section of the
total space of some disk bundle $\xi$. Let $\phi \colon \pi(S(\xi))
\rightarrow \pi(E(\xi))$ be the induced homomorphism of fundamental
groupoids. The embedding corresponds to the $LN$-situation, so we
will talk about the space $\LN_{n+1}(\phi)$.

Notice that we are provided with a map $\xi^{!} \colon \S(M)
\rightarrow \S(E(\xi),S(\xi))$. On a $k$-simplex $f \colon M'
\rightarrow M \times \Delta^k$ it is given by pulling back the disk
bundle $E(\xi) \rightarrow M$ along $f$. Let $y$ be a $0$-simplex in
$\S(E(\xi),S(\xi))$, that means $y$ is a simple homotopy equivalence
$f_y \colon (E,S) \rightarrow (E(\xi),S(\xi))$ which can also be
seen as a splitting problem along $M$. The splitting obstruction
$\theta(f_y)$ lives in the group $LN_{n}(\phi)$. Denote
\[
\S(\xi,y) = \hofibery \big(\xi^{!} \colon \S(M) \rightarrow
\S(E(\xi),S(\xi)) \big).
\]
By definition (see \cite[chapter 11]{Ran}) a $k$-simplex of
$\S(\xi,y)$ is a simple homotopy equivalence $f \colon (E,S) \times
\Delta^{k+1} \rightarrow (E(\xi), S(\xi)) \times \Delta^{k+1}$ which
respects faces, such that $f_y = \partial_\alpha f \colon (E,S)
\times \partial_\alpha \Delta^{k+1} \rightarrow (E(\xi), S(\xi))
\times \partial_\alpha \Delta^{k+1}$ for $\alpha \colon \mathbf{0}
\rightarrow \k \mathbf{+1}$, such that $\alpha(0) = k+1$, and
$\partial_{k+1} f \colon  (E,S) \times \partial_{k+1} \Delta^{k+1}
\rightarrow (E(\xi), S(\xi)) \times \partial_{k+1} \Delta^{k+1}$ is
split along $M$ where $k+1 \colon \k \rightarrow \k \mathbf{+1}$ is
obtained by omitting $k+1$ in the target.

It is not hard to see that,
\begin{itemize}
\item if $\theta(f_y) \neq 0$ then $\S(\xi,y)$ is empty,
\item if $\theta(f_y) = 0$ then $\S(\xi,y) \simeq  \LN_{n+1}(\phi)$.
\end{itemize}
For the second statement, notice that a $k$-simplex of the space
$\S(\xi,y)$ is a splitting problem of manifold $(k+3)$-ads which is
solved on the part $\partial_\alpha \cup \partial_{k+1}$ of the
boundary. As such it can be considered as a $k$-simplex of the space
$\LN_{n+1} (\phi)$, which gives an inclusion map which is the
indicated homotopy equivalence; it induces isomorphisms on homotopy
groups. See for example \cite[chapter 7.2]{Ran1} for the reference.

\subsection{The homotopy type of $F^{(1)}(V,y)$.} Consider now the
embedding $\R P(V) \hookrightarrow \R P(\VR)$. We have the
decomposition
\[
\R P(\VR) = E(\nu) \cup_{S(\nu)} C
\]
where $\nu$ is the normal disk bundle of $\R P(V)$ and $C \cong
D^n$, $n = \dim(V)$.

Consider now the embedding $\R P(V) \hookrightarrow E(\nu)$, the
induced homomorphism $\phi_V \colon \pi_1(S(\nu)) \rightarrow
\pi_1(E(\nu))$ and a homotopy equivalence $y \! = \! f_y \colon \!
(E,S) \rightarrow (E(\nu),S(\nu))$ with the splitting obstruction
$\theta (f_y) \in LN_{n-1}(\phi_V)$. Consider
\[
\S(\nu,y) = \hofibery \big(\nu^! \colon \S(\R P(V)) \rightarrow
\S(E(\nu),S(\nu))\big).
\]
By the above discussion we have that,
\begin{itemize}
\item if $\theta(f_y) \neq 0$ then $\S(\nu,y)$ is empty,
\item if $\theta(f_y) = 0$ then $\S(\nu,y) \simeq  \LN_{n}(\phi_V)$.
\end{itemize}
Now let $y = f_y \colon M \rightarrow \R P(\VR)$ be a $0$-simplex in
$\S(\R P(\VR))$. Consider the space $F^{(1)}(V,y) = \hofibery (\S(\R
P(V)) \rightarrow \S(\R P(\VR)))$. A $k$-simplex of this space can
be viewed as a splitting problem over $\R P(\VR) \times
\Delta^{k+1}$ along $\R P(V) \times \Delta^{k+1}$ satisfying certain
conditions on the boundary. The homotopy type of $F^{(1)}(V,y)$ is
identified in the following Proposition.

\begin{prop} \label{crit-F}
Let $\dim(V) = n \geq 6$ and consider the homotopy equivalence $y =
f_y \colon M \rightarrow \R P(\VR)$ as a splitting problem along $\R
P(V)$ with the splitting obstruction $\theta(f_y) \in
LN_{n-1}(\phi_V)$. Then the following holds
\begin{itemize}
\item if $\theta(f_y) \neq 0$ then $F^{(1)}(V,y)$ is empty,
\item if $\theta(f_y) = 0$ then $F^{(1)}(V,y) \simeq  \LN_{n}(\phi_V)$.
\end{itemize}
\end{prop}

\begin{proof}
Firstly notice that there is a map $\S(E(\nu),S(\nu)) \rightarrow
\S(\R P(\VR))$ given by `coning off'. This map is a homotopy
equivalence by the following argument.

By transversality there is the following homotopy cartesian square
\[
\xymatrix{ \S (\R P(\VR);E(\nu),C;S(\nu)) \ar[r] \ar[d] & \S
(C,S(\nu)) \ar[d] \\ \S (E(\nu),S(\nu)) \ar[r] & \S(S(\nu)). }
\]
By the result of \cite[Theorem 12.1]{Wall}, if $\dim(V) \geq 6$,
then the forgetful map $\S (\R P(\VR);E(\nu),C;S(\nu)) \rightarrow
\S (\R P(\VR))$ is a homotopy equivalence. It is also well known
that the forgetful map $\S (C,S(\nu)) \rightarrow \S(S(\nu))$ is a
homotopy equivalence. Therefore also the left hand vertical map is a
homotopy equivalence. The `coning off' map is its homotopy inverse.

Abusing notation we denote by $y$ also the image of $y \in \S(\R
P(\VR))$ under some homotopy inverse of the coning off map. We have
the following diagram which commutes up to homotopy
\[
 \xymatrix{
  \S(\nu,y) \ar@{-->}[d] \ar[r] & \S(\R P(V)) \ar@{=}[d]
  \ar[r]^-{\nu^{!}} & \S(E(\nu),S(\nu)) \ar[d] \\ F^{(1)}(V,y)
  \ar[r] & \S(\R P(V)) \ar[r]^-{F(\iota)} & \S(\R P(\VR)). }
\]
The left hand vertical map is induced by the right hand vertical map
which is a homotopy equivalence and so the left hand map is also a
homotopy equivalence.

By the criterion preceding the theorem we get the desired statement.
\end{proof}

\begin{prop} \label{1stder}
If $\dim(V) = n \geq 6$ we have
\[
\pi_k(\LN_n(\phi_V)) \cong LN_{n+k}(\phi_V) \cong L_{k-n}(1)
\]
\end{prop}

\begin{proof}
We have the homomorphism $\phi_V \colon 1 \rightarrow \Z_2$, where
the subscript $V$ stands for the orientation character. This equals
$(-1)^n$. By \cite[Corollary 12.9.1.]{Wall} we have that $LN_\ast (1
\rightarrow \Z_2) \cong L_{\ast + \varepsilon}(1)$ where
$\varepsilon = 0$ if the orientation character is $-1$ and
$\varepsilon = 2$ if the orientation character is $+1$. This
together with the $4$-periodicity of the $L$-groups yields the
desired statement.
\end{proof}

\subsection{An isomorphism from $\pi_k(\LN_n(\phi_V))$ to
$L_{k-n}(1)$.} Following \cite[chapter 12C]{Wall} we now give a
concrete and explicit description of the isomorphism
\[
\theta \colon \pi_k(\LN_n(\phi_V)) \rightarrow L_{k-n}(1).
\]
By definition a $k$-simplex in $\LN_n(\phi_V)$ is a splitting
problem $f \colon (E,S) \times \Delta^{k+1} \rightarrow
(E(\nu),S(\nu)) \times \Delta^{k+1}$ along $\R P(V) \times
\Delta^{k+1}$ satisfying certain conditions on the boundary. The
isomorphism $\theta$ sends the homotopy class of $f$ to the
splitting obstruction $\mathrm{rel\;} \partial \Delta^{k+1}$ of $f$
which lives in the group $L_{k-n}(1)$. Simultaneously we give a
description of the map
\[
\theta \colon \pi_0 (\S(E(\nu),S(\nu))) \rightarrow L_{-1-n}(1).
\]
A $0$-simplex of $\S(E(\nu),S(\nu))$ is a splitting problem $f
\colon (E,S) \rightarrow (E(\nu),S(\nu))$ along $\R P(V)$ and the
map $\theta$ sends $f$ to its splitting obstruction.

So we consider a splitting problem $f \colon (E,S) \times
\Delta^{k+1} \rightarrow (E(\nu),S(\nu)) \times \Delta^{k+1}$, where
$k \geq -1$ and we find a splitting obstruction $\theta(f)$ in
$L_{k-n}(1)$. If $(k-n)$ is odd, then $L_{k-n}(1)$ is trivial so it
is enough to consider the case when $(k-n)$ is even and hence also
$(n+k)$ is even.

The even dimensional $L$-groups are defined algebraically in
\cite[chapter 5]{Wall}. For a group $\pi$ with a homomorphism $w
\colon \pi \rightarrow \Z_2$ and for $2m = (k-n)$ an element of the
group $L_{k-n}(\pi)$ is represented by a {\it simple
$(-1)^m$-hermitian form} $(G,\lambda,\mu)$ where $G$ is a stably
free $\Z[\pi]$-module and the maps
\[
\lambda \colon G \times G \rightarrow \Z[\pi], \quad \mu \colon G
\rightarrow \Z[\pi]/\{x - (-1)^l \overline{x} \},
\]
where $x \mapsto \overline{x}$ is an involution on $\Z[\pi]$ induced
by the homomorphism $w$, satisfy certain conditions, see
\cite[Theorem 5.2]{Wall}. A submodule $H$ of $G$ is called a {\it
lagrangian} of the form $(G,\lambda,\mu)$ if the map $\lambda$
induces an isomorphism $G/H \rightarrow
\mathrm{Hom}_{\Z[\pi]}(H,\Z[\pi])$, $\lambda(H \times H) = 0$ and
$\mu(H)=0$. The form $(G,\lambda,\mu)$ admits a lagrangian if and
only if it represents a zero element in the $L$-group in which case
it is called a {\it hyperbolic} form.

In \cite[chapter 5]{Wall} the surgery obstruction for converting an
even dimensional degree one normal map into a simple homotopy
equivalence is described as an element of such an $L$-group.
Following \cite[chapter 12C]{Wall} we now describe how to associate
a simple $(-1)^m$-hermitian form which represents an element in
$L_{k-n}(1)$ to the splitting problem $f$ as above.

\begin{con} \label{obstruction}
So we start with $f \colon (E,S) \times \Delta^{k+1} \rightarrow
(E(\nu),S(\nu)) \times \Delta^{k+1}$ a homotopy equivalence
satisfying appropriate conditions on the boundary. Denote by $N =
f^{-1} (\R P(V) \times \Delta^{k+1})$ and by $i$ the embedding $N
\hookrightarrow (E,S) \times \Delta^{k+1}$ and note that $i$ is a
degree one normal map. The obstruction $\theta(f)$ coincides with
the obstruction to make the embedding $i$ into a homotopy
equivalence by an ambient surgery on $N$ inside $(E,S)\times
\Delta^{k+1}$ (see \cite[chapter 12C]{Wall} or \cite[chapter
7.6.]{Ran1} for the definition of the ambient surgery on a
submanifold).

This is defined as follows. Let $2l = (n+k)$. First make $i$ into an
$l$-connected map by an ambient surgery. The ordinary surgery
obstruction associated to $i$ is a simple $(-1)^l$-hermitian form
$(H, \lambda, \mu)$ over $\Z[\Z_2]$, where $H = \pi_{l+1}(i) =
\pi_{l+1}(E \times \Delta^{k+1}, N)$. As usual the elements of $H$
can be represented by $l$-dimensional spheres immersed in $N$, the
form $\lambda$ is given by the $\Z[\Z_2]$-intersection numbers of
these spheres and $\mu$ is given by the $\Z[\Z_2]$-self-intersection
numbers. Let us now instead look at the double covers and try to
change $\widetilde{i} \colon \widetilde{N} \hookrightarrow
\widetilde{E} \times \Delta^{k+1}$ by an equivariant ambient
surgery. Here $\widetilde{E}$ denotes the non-trivial double cover
of $E$. The ordinary non-equivariant surgery obstruction associated
to $\widetilde{i}$ is a simple $(-1)^l$-hermitian form $(H,
\lambda', \mu')$ over $\Z$, with the same $H$, only considered now
as a $\Z$-module. The maps $\lambda'$, $\mu'$ are now given by
$\Z$-intersection, resp. $\Z$-self-intersection numbers of the
representing immersed spheres. The form $(H,\lambda',\mu')$ is
hyperbolic. More to the point, $\widetilde{N}$ is a separating
submanifold and we have a decomposition $\widetilde{E} \! \times \!
\Delta^{k+1} = A^+ \cup_{\widetilde{N}} A^-$. The $\Z$-module $H =
\pi_{l+1}(\widetilde{i}) = \pi_{l+1}(\widetilde{E} \times
\Delta^{k+1}, \NN)$ splits as
\[
H^+ \oplus H^- = \pi_{l+1}(A^+, \widetilde{N}) \oplus \pi_{l+1}(A^-,
\widetilde{N}).
\]
The involution $T$ switches $H^+$ and $H^-$ which are both
lagrangians of $(H, \lambda', \mu')$.

The equivariant ambient surgery obstruction for $\widetilde{i}$ is
now defined as follows. Denote by $\lambda_0 (x,y) = \lambda'
(x,ty)$, $\mu_0(x) = t\mu (x)$, where $t$ is the nontrivial element
of $\Z_2$. Then $(H^+, \lambda_0, \mu_0)$ is a simple
$(-1)^{m}$-hermitian form over $\Z$ where $2m = (k-n)$. It
represents the desired obstruction $\theta(f)$ in $L_{k-n}(1)$.

Furthermore we may assume that $H^+ \cong \pi_{l+1} (A^+,\NN)$ is a
free $\Z$-module on generators $\{ e_i \}_{i \in \r}$ whose
boundaries are represented by the spheres $S^l_i$ disjointly
embedded in $\NN$ (because $H^+$ is a lagrangian over $\Z$). Let
$\NN \times I$ be the normal disk bundle of $\NN$ in $\widetilde{E}
\times \Delta^{k+1}$. Then the generators $e_i$ can be represented
by the disks $D^{l+1}_i$ disjointly embedded in $\widetilde{E}
\times \Delta^{k+1}$ so that the disk $D^{l+1}_i$ is attached to the
normal bundle $\NN \times I$ at the sphere $S^l_i \times \{1\}$ (see
\cite{BL}). Similarly $H^-$ is generated by the elements $\{t
e_i\}_{i \in \r}$ whose boundaries are represented by the spheres $T
S^l_i$ and the elements $t e_i$ are themselves represented by the
disks $T D^{l+1}_i$ attached to $\NN \times I$ at $T S^l_i \times
\{-1\}$.
\end{con}

\begin{rem}
We did not prove that $\theta$ is an isomorphism if $k \geq 0$, we
have only given a description of the map $\theta$. However, from the
fact that $\theta$ is an isomorphism if $k \geq 0$ follows an
important realization theorem which says that if $k \geq 0$ then any
element of $L_{k-n}(1)$ can be represented by a splitting problem
over $(E(\nu),S(\nu)) \times \Delta^{k+1}$ along $\R P(V) \times
\Delta^{k+1}$. We refer to \cite[chapter 12C]{Wall} for the details.
\end{rem}

\begin{rem} \label{equiv-model-of-F1}
In this section we have used a model of $F^{(1)}(V,y)$ such that a
$k$-simplex was a splitting problem over $\R P(\VR) \times
\Delta^{k+1}$ along $\R P(V) \times \Delta^{k+1}$ satisfying certain
conditions on the boundary. In the next section it will be more
convenient to use another model where a $k$-simplex is an
equivariant splitting problem over $S(\VR) \times \Delta^{k+1}$
along $S(V) \times \Delta^{k+1}$ (it is clear what involutions are
meant) satisfying corresponding conditions on the boundary.
\end{rem}

%%%%%%%%%%%%%%%%%%%%%%%%%%%%%%%%%%%%%%%%%%%%%%%%%%%%%%%%%%%%%%%%

\section{The Functor $F$ is polynomial of degree $\leq 1$} \label{chap:ThetaF1}

In this section we prove the two main theorems of this paper,
Theorem \ref{thm1}. and Theorem \ref{thm2}. We start with:

\begin{thm} \label{thm1}
Let $F \colon \sJ \rightarrow \Spaces_\ast$ be the functor defined
by $V \mapsto \S(\R P(V))$. If $\dim(V) \geq 6$ then the canonical
maps
\begin{equation} \label{map}
F(V) \rightarrow \holimsub{0 \neq U \subset \R^2} F(\VU)
\end{equation}
are homotopy equivalences.
\end{thm}

An immediate corollary of Theorem \ref{thm1} is the following.
\begin{cor} \label{cor1}
If $\dim(V) \geq 6$ then the canonical maps
\[
F(V) \rightarrow T_1F(V)
\]
are homotopy equivalences.
\end{cor}

\begin{proof}[Proof of Corollary \ref{cor1}]
The definition of the functor $T_1E$ for a general continuous
functor $E \colon \sJ \rightarrow \Spaces_\ast$ is given in
\cite[Theorem 6.3]{Weiss}. It says that $T_1E$ is the homotopy
colimit of the direct system
\[
E \xrightarrow{\rho} \tau_1 E \xrightarrow{\tau_1(\rho)} (\tau_1)^2
E \xrightarrow{(\tau_1)^2(\rho)} \cdots,
\]
where
\[
\tau_1 E (V) = \holimsub{0 \neq U \subset \R^2} E(\VU),
\]
and the natural transformation $\rho \colon E \rightarrow \tau_1 E$
is given by the canonical maps. Theorem \ref{thm1} says that in the
case of the functor $F$ the maps $\rho \colon F(V) \rightarrow
\tau_1 F(V)$ are homotopy equivalences if $\dim (V) \geq 6$. The
claim of Corollary \ref{cor1} follows.
\end{proof}

\noindent \textbf{Overview of the proof of Theorem \ref{thm1}.} The
details of the proof of Theorem \ref{thm1} are rather complicated.
Therefore we first give an overview of the proof.

The pattern of the proof is to verify the condition of Theorem
\ref{poly1} for the map (\ref{map}) to be a homotopy equivalence, if
$\dim(V) \geq 6$. This condition says that we have to verify that a
certain map $\ssigma_1 \colon F^{(1)}(V,y) \rightarrow X_y$ is a
homotopy equivalence for all choices of $y \in F(\VR)$.

The space $X_y$ and the map $\ssigma_1$ were described for a general
continuous functor $E \colon \sJ \rightarrow \Spaces_\ast$ in
Construction \ref{map2}. We work with the $\Delta$-set models, so in
Construction \ref{map3} we provide a $\Delta$-set model of the space
$X_y$ and we also provide a $\Delta$-set description of the map
$\ssigma_1$ which up to homotopy commutes with the geometric
realization.

The homotopy type of the space $F^{(1)}(V,y)$ has been described in
Proposition \ref{crit-F}. It says that $F^{(1)}(V,y)$ is empty if
and only if the splitting obstruction $\theta(f_y)$ of a certain
equivariant splitting problem $f_y$ over $S(\VR)$ along $S(V)$ is
non-zero (see also Remark \ref{equiv-model-of-F1}). If $\theta(f_y)
= 0$ then the space $F^{(1)}(V,y)$ is homotopy equivalent to a
certain space $\LN_n(\phi_V)$ whose homotopy groups are calculated
in Proposition \ref{1stder}.

We give a corresponding description of the homotopy type of the
space $X_y$ in Proposition \ref{crit-X}. It says that $X_y$ is empty
if and only if the splitting obstruction $\theta(f'_y)$ of a certain
equivariant splitting problem $f'_y$ over $S(\VR^2) \times I$ and
along $S(\VR) \times I$ is non-zero. If $\theta(f'_y) = 0$ then the
space $X_y$ is homotopy equivalent to the space $\Omega
\LN_{n+1}(\phi_{\VR})$ whose homotopy groups are also calculated by
Proposition \ref{1stder}.

As we have indicated in Remark \ref{pattern} the proof then splits
into two parts. To prove the condition (1) we have to show that the
space $F^{(1)}(V,y)$ is empty if and only if the space $X_y$ is
empty. This is the content of Proposition \ref{prop1}. It is proved
using the criteria mentioned above. Notice that the equivariant
splitting obstructions $\theta(f_y)$ and $\theta(f'_y)$ live in the
same group $L_{-1-n}(1)$. Therefore it is enough to check that
$\theta(f_y) = \theta(f'_y)$ in $L_{-1-n}(1)$. The proof of this
statement is the crucial part of the whole proof. Unfortunately it
is rather technical. It is given in a sequence of lemmas: Lemma
\ref{lem1}, \ref{lem2}, \ref{lem3}.

To prove the condition (2) we have to show that if $F^{(1)}(V,y)$ is
not empty then the map $\ssigma_1 \colon F^{(1)}(V,y) \rightarrow
X_y$ is a homotopy equivalence. This is done by showing that it
induces isomorphisms on homotopy groups. Using Propositions
\ref{crit-F}, \ref{crit-X} and \ref{1stder} we see that the spaces
$F^{(1)}(V,y)$ and $X_y$ have isomorphic homotopy groups. Thus we
have to show that the map $\ssigma_1$ induces isomorphisms between
abstractly isomorphic groups. The idea here is the following. An
element of $\pi_k(F^{(1)}(V,y)) \cong L_{k-n}(1)$ can be represented
by a splitting obstruction $\theta(f)$ of some equivariant splitting
problem $f$ over $S(\VR)\times \Delta^{k+1}$ along $S(V) \times
\Delta^{k+1}$. The image $\ssigma_1(f)$ is a splitting problem over
$S(\VR^2) \times \Delta^{k+1} \times I$ along $S(\VR) \times
\Delta^{k+1} \times I$ with the splitting obstruction
$\theta(\ssigma_1(f))$ which represents the homotopy class of
$\ssigma_1(f)$ in $\pi_k(X_y) \cong L_{k-n}(1)$. The condition (2)
is proved by showing that $\theta(f) = \theta(\ssigma_1(f))$ in
$L_{k-n}(1)$. This statement is the content of Proposition
\ref{prop2}. Its proof turns out to be just a parameterized version
of the proof of the statement of condition (1).

\begin{con} \label{map3}
We need a $\Delta-\mathrm{set}$ model of the space $X_y$ and a
description of the map $\ssigma_1$ on the level of
$\Delta-\mathrm{sets}$. The descriptions given here and the
descriptions from Construction \ref{map2} commute with the geometric
realization up to homotopy.

Recall the $\Delta-\mathrm{set}$ model for $F^{(1)}(V,y)$. Here $y =
f_y \colon S(\VR) \rightarrow S(\VR)$ is an equivariant splitting
problem along $S(V)$. A $k$-simplex of $F^{(1)}(V,y)$ is an
equivariant homotopy equivalence $f \colon S(\VR) \times
\Delta^{k+1} \rightarrow S(\VR) \times \Delta^{k+1}$ w.r.t. the
antipodal involution on the target and w.r.t. some free involution
on the source, it respects faces, it is transverse to $S(V) \times
\Delta^{k+1}$, for $\alpha \colon \mathbf{0} \rightarrow \k
\mathbf{+1}$, $\alpha(0)=k+1$, we have $f_y =
\partial_\alpha f \colon S(\VR) \times \partial_\alpha \Delta^{k+1}
\rightarrow S(\VR) \times \partial_\alpha \Delta^{k+1}$ and
$\partial_{k+1} f$ is equivariantly split along $S(V) \times
\Delta^k$.

Now we give a $\Delta-\mathrm{set}$ model for $X_y$. A $k$-simplex
of $X_y$ is an equivariant homotopy equivalence $f' \colon S(\VR^2)
\times \Delta^{k+1} \times I \rightarrow S(\VR^2) \times
\Delta^{k+1} \times I$ w.r.t. the antipodal involution on the target
and w.r.t. to some free involution on the source which respects
faces and satisfies the following conditions on the boundary.

Firstly define
\[
\partial_0 (\Delta^{k+1} \times I) = (\partial_{k+1} \Delta^{k+1}
\times I) \cup (\Delta^{k+1} \times \partial I),
\]
\[
\partial_1 (\Delta^{k+1} \times I) = \partial_\alpha \Delta^{k+1} \times I,
\]
where $\alpha \colon \mathbf{0} \rightarrow \mathbf{k+1}$ is again
given by $\alpha(0) = k+1$.

The conditions on $f'$ are the following. Firstly $f'$ restricted to
$S(\VR^2) \times \partial_0 (\Delta^{k+1} \times I)$ is split along
$S(\VR) \times \partial_0 (\Delta^{k+1} \times I)$. The restriction
of $f'$ to $S(\VR^2) \times \partial_1 (\Delta^{k+1} \times I)$ is
fixed and it is given as the composition
\begin{equation} \label{map-f'_y}
f'_y = (S(\zeta))^{-1} \circ (\Sigma(f_y) \times \id) \circ
S(\zeta).
\end{equation}
Here $\Sigma$ means the unreduced suspension map and the map
$S(\zeta) \colon S(\VR^2) \times I \rightarrow S(\VR^2) \times I$ is
a homeomorphism given by $(x,t) \mapsto (\zeta_t(x),t)$ where
$\zeta_t$ is the unique rotation of $\VR^2$ which fixes $V$ and
sends $(0,1,0)$ to $(0,t,\sqrt{1-t^2})$.

The map $\ssigma_1 \colon F^{(1)}(V,y) \rightarrow X_y$ is given by
the rule
\begin{equation} \label{map-f'}
\ssigma_1 \colon f \mapsto (S(\zeta))^{-1} \circ (\Sigma(f) \times
\id) \circ S(\zeta).
\end{equation}

\begin{figure}[!htp]
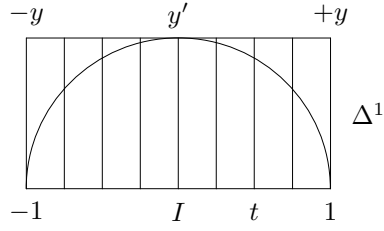

\[
\xy 0;<2cm,0cm>: (0,0)*\cir<2cm>{u^d}, (-1,0) ; (1,0) **@{-}, (-1,1)
; (1,1) **@{-}, (-1,0) ; (-1,1) **@{-}, (1,0) ; (1,1) **@{-}, (-0.75,0) ;
(-0.75,1) **@{-}, (-0.5,0) ; (-0.5,1) **@{-}, (-0.25,0) ; (-0.25,1)
**@{-}, (0,0) ; (0,1) **@{-}, (0.25,0) ; (0.25,1) **@{-}, (0.5,0) ;
(0.5,1) **@{-}, (0.75,0) ; (0.75,1) **@{-},
,(0,-0.15)*={I};(1.25,0.5)*={\Delta^1};
(-1,-0.15)*={-1};(0.5,-0.15)*={t}; (1,-0.15)*={1};
(-1,1.15)*={-y};(0,1.15)*={y'}; (1,1.15)*={+y}
\endxy
\]
\caption{Schematic picture of a $0$-simplex in $X_y$ and of the map
$\ssigma_1$.} \label{fig:sigma3}
\end{figure}
A $0$-simplex of $X_y$ is depicted in Figure \ref{fig:sigma3}. The
top edge corresponds to the part $\partial_1(\Delta^1 \times I)$ of
the boundary. The closure of the complement of $\partial_1(\Delta^1
\times I)$ in the boundary is $\partial_0 (\Delta^1 \times I)$. On
each vertical line over $t \in I$ the restriction of the map
$\ssigma_1(f)$ is given by $(\zeta_t)^{-1} \circ \Sigma(f) \circ
\zeta_t$. Compare this with Figure \ref{fig:sigma2}. The vertical
lines here correspond after reparametrization to the rays from the
origin in Figure \ref{fig:sigma2}.
\end{con}

Note that the map $f'_y$ can be seen as an equivariant splitting
problem along $S(\VR) \times I$ with the splitting obstruction
$\theta(f'_y) \in LN_{n+1}(\phi_{\VR}) \cong L_{-1-n}(1)$. We need
the following criterion.

\begin{prop} \label{crit-X}
We have that
\begin{itemize}
\item if $\theta(f'_y) \neq 0$ then the space $X_y$ is empty,
\item if $\theta(f'_y) = 0$ then $X_y \simeq \Omega \LN_{n+1}(\phi_{\VR})$.
\end{itemize}
\end{prop}

\begin{proof}
If $X_y$ is not empty then there exists a $0$-simplex, say $f$, of
$X_y$. Since $f$ restricted to $(S(\VR^2) \times \partial_0
(\Delta^{k+1} \times I))$ is by definition split along $S(\VR)
\times \partial_0 (\Delta^{k+1} \times I)$ any such a $0$-simplex
$f$ can be seen as a solution of the splitting problem $f'_y$ and
hence $\theta(f'_y) = 0$. As we have remarked in Construction
\ref{map2} if the space $X_y$ is not empty it is homotopy equivalent
to a certain space of paths in the space $F^{(1)}(\VR,y') \simeq
\LN_{n+1}(\phi_{\VR})$. It is not difficult to see that this path
space is homotopy equivalent to $\Omega \LN_{n+1}(\phi_{\VR})$.
\end{proof}

Now we are ready to proceed with a verification of the condition of
Theorem \ref{poly1} for the functor $F$ to be polynomial of degree
$\leq 1$. We start with the condition (1).

\begin{prop} \label{prop1}
$F^{(1)}(V,y)$ is empty if and only if $X_y$ is empty.
\end{prop}

Propositions \ref{crit-F}, \ref{crit-X} tell us when the spaces
$F^{(1)}(V,y)$ and $X_y$ are empty. We see that both obstructions
$\theta(f_y)$ and $\theta(f'_y)$ live in the group $L_{-1-n}(1)$. In
order to prove the condition (1) it is enough to show that
$\theta(f_y) = \theta (f'_y)$ in $L_{-1-n}(1)$. This is done in the
sequence of the following three lemmas.

\begin{rem}
Note that a precomposition with a homeomorphism does not change the
splitting obstruction of a splitting problem. Therefore when
investigating $\theta(f'_y)$ it is enough to look at the splitting
obstruction of the map $(S(\zeta))^{-1} \circ (\Sigma (f_y) \times
\id)$ and so from now on we denote
\begin{equation} \label{mapf'_y2}
f'_y = (S(\zeta))^{-1} \circ (\Sigma (f_y) \times \id).
\end{equation}
\end{rem}

\begin{lem} \label{lem1}
The map $f'_y \colon S(\VR^2) \times I \rightarrow S(\VR^2) \times
I$ is transverse to the submanifold $S(\VR) \times I$.
\end{lem}

\begin{proof}
One would expect that in order to find the splitting obstruction
$\theta(f'_y)$ it is necessary to first adjust the map $f'_y$ to
make it transverse to $S(\VR) \times I$. Amazingly, the map $f'_y$
as constructed by the formula (\ref{mapf'_y2}) already is transverse
to $S(\VR) \times I$! This is what we prove in this lemma, more
precisely we show an equivalent statement that the map $(\Sigma
(f_y) \times \id)$ is transverse to $S(\zeta) \circ j$ where $j$
denotes the embedding $S(\VR) \times I \hookrightarrow S(\VR^2)
\times I$.

We use an easier notion of transversality which is good enough for
our purposes. It goes as follows. Let $M \hookrightarrow X$ be a
codimension $1$ embedding. Suppose that there exist an open
neighborhood $E$ of $M$ in $X$ and a submersion $g \colon E
\rightarrow \R$ such that $g^{-1}(0) = M$. (This implies that $M$
has a trivial normal bundle.) We say that a map $f \colon Y
\rightarrow X$ is transverse to $M$ if $g \circ f \colon f^{-1}(E)
\rightarrow \R$ is a submersion in a neighborhood of $f^{-1}(M)$.
Note that this condition can be verified locally.

In our case, $X = S(\VR^2) \times I$ and $M$ is the image of $S(\VR)
\times I$ under the embedding $S(\zeta) \circ j$. Let $g = h \circ
(S(\zeta))^{-1}$ where $h \colon S(\VR^2) \times I \rightarrow \R$
is the `height function' defined by the projection on the second
coordinate of $\R^2$. (This should be suitably restricted to a
neighborhood of $M$ in $X$.)

We use the following decompositions (\ref{decomposition1}),
(\ref{decomposition2}) into (overlapping) codimension $0$
(non-compact) submanifolds
\begin{equation} \label{decomposition1}
S(\VR) = E(\widetilde{\nu}) \cup C_\pm
\end{equation}
where $E(\widetilde{\nu})$ is the total space of the open normal
vector bundle $\widetilde{\nu}$ of $S(V)$ in $S(\VR)$, and $C_\pm =
D_\pm(\VR) \smallsetminus S(V)$, where $D_\pm(\VR) = S(\VR) \cap
(\VR_\pm)$ with $\R_\pm = \{ u \in \R \, | \, \pm u \geq 0 \}$, and
\begin{equation} \label{decomposition2}
S(\VR^2) = E(\widetilde{\mu}) \cup D_\pm
\end{equation}
where $E(\widetilde{\mu})$ is the total space of the open normal
vector bundle $\widetilde{\mu}$ of $S(V)$ in $S(\VR^2)$ and  $D_\pm
= D_\pm(\VR^2) \smallsetminus S(V)$, where $D_\pm(\VR^2) = S(\VR^2)
\cap (\VR \! \oplus \! \R_\pm)$.

Notice that the embedding $S(\zeta) \circ j$ restricts to the
embeddings
\[
E(\widetilde{\nu}) \times I  \hookrightarrow E(\widetilde{\mu})
\times I, \quad C_\pm \times I \hookrightarrow D_\pm \times I.
\]

We give two different arguments for the transversality on these two
kinds of codimension $0$ submanifolds of $S(\VR^2) \times I$

We start with the embedding $E(\widetilde{\nu}) \times I
\hookrightarrow E(\widetilde{\mu}) \times I$. This is an embedding
of a codimension $1$ subbundle where both bundles are over $S(V)
\times I$. The map $g$ restricted to $E(\widetilde{\mu}) \times I$
can be described as $g \colon (x,t) \mapsto (x - \pr_l(x))$ where
$\pr_l$ is the orthogonal projection in the fiber of
$E(\widetilde{\mu}) \times \{t\}$ onto the image of
$E(\widetilde{\nu}) \times \{t\}$.

Now the map $f_y$ was transverse to $S(V)$, therefore the map
$\Sigma (f_y)$ is also transverse to $S(V)$. Let $\NN :=
(f_y)^{-1}(S(V))$ and let $\widetilde{\nu}_{\NN}$ be the open normal
bundle of $\NN$ in $S(\VR)$ and denote by $\widetilde{\mu}_{\NN}$
the open normal bundle of $\NN$ in $S(\VR^2)$. The restriction
$\Sigma (f_y)| \colon E(\widetilde{\mu}_{\NN}) \rightarrow
E(\widetilde{\mu})$ is a bundle map. Now it is not difficult to see
that the restriction $\Sigma (f_y)| \times \id \colon
E(\widetilde{\mu}_{\NN}) \times I \rightarrow E(\widetilde{\mu})
\times I$ is transverse to $E(\widetilde{\nu}) \times I$ in the
sense above.

Next we look at the embedding $C_\pm \times I \hookrightarrow D_\pm
\times I$. Observe that the composition $C_\pm \times I
\hookrightarrow D_\pm \times I \rightarrow D_\pm$ is a
homeomorphism. Therefore we can consider the embedding $C_\pm \times
I \cong D_\pm \hookrightarrow D_\pm \times I$ as an embedding of a
graph of some function $\omega \colon D_\pm \rightarrow I$ into
$D_\pm \times I$. The restriction of $g$ can be expressed as $(x,t)
\mapsto t - \omega(x)$. Again it is not difficult to see that the
restriction $\Sigma (f_y)| \times \id \colon (D_\pm(\VR^2)
\smallsetminus \NN) \times I \rightarrow D_\pm \times I$ is
transverse to $C_\pm \times I$ in the sense above.

The details are left to the reader.
\end{proof}

In the next lemma we identify the homotopy type of $\NN' :=
(f'_y)^{-1} (S(\VR) \times I)$.

\begin{lem} \label{lem2}
Let $\tau \colon \NN' \rightarrow \R$ be the composition $\NN'
\hookrightarrow S(\VR^2) \times I \rightarrow I \hookrightarrow \R$.
Then $\tau$ is transverse to $s \in \R$ for all $s \neq 0$. Hence
$\NN' \simeq \tau^{-1}(0) \cong \Sigma \NN$ where $\Sigma \NN$ means
the unreduced suspension of $\NN$.
\end{lem}

\begin{proof}
A key observation here is that the manifold (with boundary) $\NN'$
has a decomposition
\begin{equation} \label{decomposition3}
\NN' = E(\widetilde{\nu}_{\NN}) \times I \cup (D_\pm(\VR^2)
\smallsetminus \NN)
\end{equation}
into codimension $0$ submanifolds which corresponds to the previous
decomposition $S(\VR) \times I = (E(\widetilde{\nu}) \times I) \cup
D_\pm$. This can be seen as follows.

We can think of $\NN'$ as a pullback of $S(\zeta) \circ j$ along
$\Sigma (f_y) \times \id$:
\[
\xymatrix@!C=8.8em{ \NN' \ar[d] \ar[r] & S(\VR) \times I \
\ar[d]^{S(\zeta) \circ j} \\ S(\VR^2) \times I \ar[r]^{\Sigma (f_y)
\times \id} & S(\VR^2) \times I.}
\]
Look first at the following restriction of this diagram:
\[
\xymatrix@!C=8.8em{ ? \ar[d] \ar[r] & E(\widetilde{\nu}) \times I \
\ar[d]^{S(\zeta) \circ j} \\ E(\widetilde{\mu}_{\NN}) \times I
\ar[r]^{\Sigma (f_y)| \times \id} & E(\widetilde{\mu}) \times I.}
\]
In this diagram the lower horizontal arrow is a bundle map over the
restriction $f_y | \times \id \colon \NN \times I \rightarrow S(V)
\times I$. Here $E(\widetilde{\mu}_{\NN}) \times I$ is a
$2$-dimensional vector bundle over $\NN \times I$,
$E(\widetilde{\mu}) \times I$ is a $2$-dimensional vector bundle
over $S(V) \times I$ and $E(\widetilde{\nu}) \times I$ is its
codimension $1$ subbundle. It follows that the top row $?
\rightarrow E(\widetilde{\nu}) \times I$ is a bundle map and it also
follows from the construction that $?$ can be identified with
$E(\widetilde{\nu}_{\NN}) \times I$.

Secondly, look at the restriction of the diagram to $C_\pm \times I
\cong D_\pm$. Here, because the map $S(\zeta) \circ j$ is a map over
$I$, we can look at the following pullback diagrams
\[
\xymatrix@!C=8.8em{ ?? \ar[d] \ar[r] & C_\pm \times I \
\ar[d]^{\pr_1 \circ S(\zeta) \circ j} \\ S(\VR^2) \ar[r]^{\Sigma
(f_y)} & S(\VR^2).}
\]
where $\pr_1$ is the projection onto the first factor. As noted
before the composition $(\pr_1 \circ S(\zeta) \circ j) \colon C_\pm
\times I \rightarrow D_\pm$ is a homeomorphism and thus we obtain
that $?? \cong \Sigma (f_y)^{-1}(D_\pm) = (D_\pm(\VR^2)
\smallsetminus \NN)$. This yields the decomposition
(\ref{decomposition3}).

We describe the map $\tau$ on each submanifold in the decomposition.
The restriction of  $\tau$ to the submanifold
$E(\widetilde{\nu}_{\NN}) \times I$ is just the projection:
\[
\tau| \; \colon E(\widetilde{\nu}_{\NN}) \times I \rightarrow I
\hookrightarrow \R
\]
which is clearly a submersion.

On the two submanifolds $D_\pm(\VR^2) \smallsetminus \NN$ the
restriction of the map $\tau$ can be described as follows. Observe
first that $D_\pm(\VR^2) \smallsetminus \NN \subset D_\pm(\VR^2)$,
and $D_\pm \subset D_\pm(\VR^2)$ and that $D_\pm(\VR^2)$ can be seen
as a cone on $S(\VR)$. The restriction of $\tau$ is the following
composition:
\[
\tau| \; \colon  D_\pm(\VR^2) \smallsetminus \NN
\xrightarrow{C(f_y)} D_\pm \cong C_\pm \times I \rightarrow I
\hookrightarrow \R
\]
where $C(f_y)$ denotes the canonical extension of $f_y$ to the cones
suitably restricted. The homeomorphism $D_\pm \cong C_\pm \times I$
can be described as follows. The fiber $C_\pm \times \{s\}$ is
embedded in $D_\pm$ as an open `meridian' through the point
$(0,s,\pm\sqrt{1-s^2})$. It is not difficult to see that the cone
map $C(f_y)$ is transverse to all such meridians except the one
corresponding to $s=0$ (details left to the reader).

The identification $\tau^{-1}(0) \cong \Sigma \NN$ follows from the
construction of the map $f'_y$ by formula (\ref{mapf'_y2}).
\end{proof}

\begin{lem} \label{lem3}
Let $2m = -1-n$ and let $\theta(f_y) \in L_{-1-n}(1)$ be represented
by a simple $(-1)^m$-hermitian form $(H^+,\lambda_0,\mu_0)$. Then
$\theta(f'_y) \in L_{-1-n}(1)$ is represented by the same simple
$(-1)^m$-hermitian form.
\end{lem}

\begin{figure}[!htp]
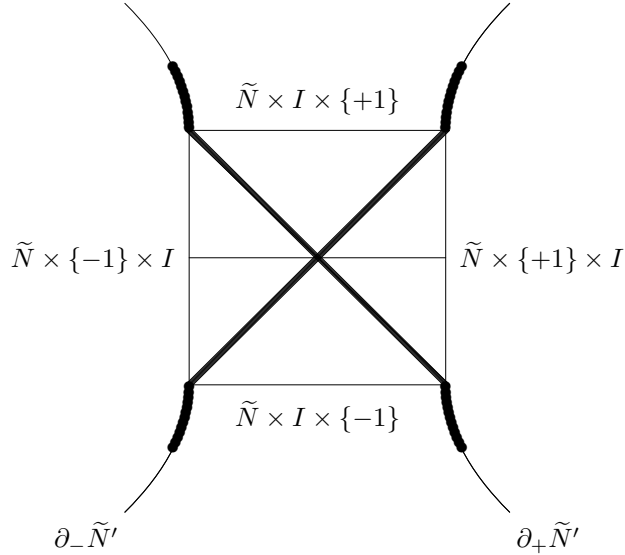

\[
\xy/r4pc/:,(-1.5,-2)="A0", (1.5,-2)="A1", (1.5,2)="A2",
(-1.5,2)="A3", (-1,-1)="B0", (1,-1)="B1", (1,1)="B2", (-1,1)="B3",
(-1,0)="C0", (0,0)="C1", (1,0)="C2", (-0.99,-1)="B00",
(1.01,1)="B20", (-0.98,-1)="B01", (1.02,1)="B21", (-0.97,-1)="B02",
(1.03,1)="B22", (-1.01,-1)="B03", (0.99,1)="B23", (-1.02,-1)="B04",
(0.98,1)="B24", (-1.03,-1)="B05", (0.97,1)="B25", (1.03,-1)="B10",
(-0.97,1)="B30", (1.02,-1)="B11", (-0.98,1)="B31", (1.01,-1)="B12",
(-0.99,1)="B32", (0.99,-1)="B13", (-1.01,1)="B33", (0.98,-1)="B14",
(-1.02,1)="B34", (1.03,-1)="B15", (-1.03,1)="B35", "B0" ; "B1"
**@{-}, "B1" ; "B2"
**@{-}, "B2" ; "B3" **@{-}, "B0" ; "B3" **@{-}, "C0" ; "C2"
**@{-}, "B0" ; "B2" **@{-}, "B00" ; "B20" **@{-}, "B01" ; "B21"
**@{-}, "B02" ; "B22" **@{-}, "B03" ; "B23" **@{-}, "B04" ; "B24"
**@{-}, "B05" ; "B25" **@{-},  "B1" ; "B3" **@{-}, "B10" ; "B30"
**@{-}, "B11" ; "B31" **@{-}, "B12" ; "B32" **@{-}, "B13" ; "B33"
**@{-}, "B14" ; "B34" **@{-}, "B15" ; "B35" **@{-}, "A0" ; "B0"
**\crv{(-1,-1.5)} ?(0.5)="D0", "A1" ; "B1" **\crv{(1,-1.5)}
?(0.5)="D1", "A2" ; "B2" **\crv{(1,1.5)} ?(0.5)="D2", "B3" ; "A3"
**\crv{(-1,1.5)} ?(0.5)="D3" *{\bullet}, ?(0.46)*{\bullet},
?(0.42)*{\bullet}, ?(0.38)*{\bullet}, ?(0.34)*{\bullet},
?(0.30)*{\bullet}, ?(0.26)*{\bullet}, ?(0.22)*{\bullet},
?(0.18)*{\bullet}, ?(0.14)*{\bullet}, ?(0.1)*{\bullet},
?(0.06)*{\bullet}, ?(0.02)*{\bullet}, "B2" ; "A2" **\crv{(1,1.5)}
?(0.5)="D2" *{\bullet}, ?(0.46)*{\bullet}, ?(0.42)*{\bullet},
?(0.38)*{\bullet}, ?(0.34)*{\bullet}, ?(0.30)*{\bullet},
?(0.26)*{\bullet}, ?(0.22)*{\bullet}, ?(0.18)*{\bullet},
?(0.14)*{\bullet}, ?(0.1)*{\bullet}, ?(0.06)*{\bullet},
?(0.02)*{\bullet}, "B1" ; "A1" **\crv{(1,-1.5)} ?(0.5)="D1"
*{\bullet}, ?(0.46)*{\bullet}, ?(0.42)*{\bullet},
?(0.38)*{\bullet}, ?(0.34)*{\bullet}, ?(0.30)*{\bullet},
?(0.26)*{\bullet}, ?(0.22)*{\bullet}, ?(0.18)*{\bullet},
?(0.14)*{\bullet}, ?(0.1)*{\bullet}, ?(0.06)*{\bullet},
?(0.02)*{\bullet}, "B0" ; "A0" **\crv{(-1,-1.5)} ?(0.5)="D0"
*{\bullet}, ?(0.46)*{\bullet}, ?(0.42)*{\bullet},
?(0.38)*{\bullet}, ?(0.34)*{\bullet}, ?(0.30)*{\bullet},
?(0.26)*{\bullet}, ?(0.22)*{\bullet}, ?(0.18)*{\bullet},
?(0.14)*{\bullet}, ?(0.1)*{\bullet}, ?(0.06)*{\bullet},
?(0.02)*{\bullet}, (-1.8,-2.2)*=0{\partial_{-}\widetilde{N}'},
(1.8,-2.2)*=0{\partial_{+}\widetilde{N}'}, (1.75,0)*=0{\widetilde{N}
\times \{+1\} \times I}, (-1.75,0)*=0{\widetilde{N} \times \{-1\}
\times I}, (0,-1.25)*={\widetilde{N} \times I \times \{-1\}},
(0,1.25)*={\widetilde{N} \times I \times \{+1\}}
\endxy
\]
\caption{Intersections of the generators of $(H, \lambda, \mu)$ in
$\NN'$.} \label{intersections}
\end{figure}

\begin{proof}
Recall that $f_y \colon S(\VR) \rightarrow S(\VR)$ is a homotopy
equivalence about which we can assume that it is transverse to
$S(V)$ and also that it is obtained from a homotopy equivalence
$E(\widetilde{\nu}) \rightarrow E(\widetilde{\nu})$ by `coning off'
(by Proposition \ref{crit-F}). Here $E(\widetilde{\nu})$ denotes the
total space of the normal disk bundle of $S(V)$ in $S(\VR)$. The
splitting obstruction $\theta(f_y)$ is represented by a simple
$(-1)^m$-hermitian form $(H^+,\lambda_0,\mu_0)$ obtained as an
equivariant ambient surgery obstruction associated to the embedding
$\widetilde{i} \colon \NN = (f_y)^{-1}(S(V)) \hookrightarrow
E(\widetilde{\nu})$ (see Construction \ref{obstruction}).

The homotopy equivalence $f'_y \colon S(\VR^2) \times I \rightarrow
S(\VR^2) \times I$ is constructed by formula (\ref{mapf'_y2}). Let
$E(\widetilde{\nu}')$ be the total space of the normal disk bundle
of $S(\VR)$ in $S(\VR^2)$ and note that $f'_y$ already is in the
form of a homotopy equivalence obtained from another homotopy
equivalence $E(\widetilde{\nu}') \times I \rightarrow
E(\widetilde{\nu}') \times I$ by `coning off' so there is no need to
adjust it. The splitting obstruction $\theta(f'_y)$ should now be
read off as the equivariant ambient surgery obstruction from the
embedding $\widetilde{i}' \colon \NN' \hookrightarrow
E(\widetilde{\nu}') \times I$ by the recipe given in Construction
\ref{obstruction}.

In previous Lemma \ref{lem2} we have made a first step by
identifying the homotopy type of $\NN'$. Moreover it follows from
the construction of the map $f'_y$ by formula (\ref{mapf'_y2}) that
the embedding $\widetilde{i}'$ when restricted to $\tau^{-1}(0)$ is
just the unreduced suspension map $\Sigma \widetilde{i}$. So we have
a commutative diagram
\[
\xymatrix{ \Sigma \NN \ar[d]^{\simeq} \ar[r]^-{\Sigma \widetilde{i}}
& E(\widetilde{\nu}') \times \{0\} \ \ar[d]^{\simeq} \\ \NN'
\ar[r]^-{\widetilde{i}'} &  E(\widetilde{\nu}') \times I}
\]
where all maps are embeddings.

Let $2l = n-1$. From the diagram above it follows that
$\widetilde{i}'$ is $(l+1)$-connected, and also that
$\pi_{l+2}(\widetilde{i}') \cong \pi_{l+1}(\widetilde{i}) \cong H$
and that the corresponding splitting of $H$ is the same $H^+ \oplus
H^-$ as before. It only remains to find the embedded spheres
representing the generators of $H$ and to show that their
intersection numbers are preserved.

Let $\{e_i, te_i\}$ be the set of generators of $H$ as a
$\Z$-module, where generators are represented by embedded spheres
$S_i, TS_i \hookrightarrow \NN$ as in Construction
\ref{obstruction}. Denote the images under suspension map by $S'_i,
TS'_i \hookrightarrow \Sigma \NN \hookrightarrow \NN'$. We need to
put them in a general position.

First we describe the topology of $\NN'$ more. The manifold $\NN'$
is a special manifold $3$-ad with two distinguished parts of the
boundary $\partial_\pm \NN' = S(\VR) \times \{\pm 1\}$. From the
decomposition (\ref{decomposition3}) we also see that the product
$\NN \times I$ is a separating submanifold of $\NN'$. By taking its
normal disk bundle (using collars) we obtain the codimension $0$
embedding $\NN \times I^2 \hookrightarrow \NN'$, which when
restricted to $\NN \times \{\pm 1\} \times I \hookrightarrow
\partial_\pm \NN'$ coincides with the embedding of the closed normal bundle
of $\NN \times \{\pm 1\} \hookrightarrow \partial_\pm \NN'$. In
particular there are disks $D^\pm_i, TD^\pm_i$ attached to the
normal bundle $\NN \times I$ as described in Construction
\ref{obstruction}. See Figure \ref{intersections}.

Denote now by $d_\pm \colon I \rightarrow I^2$ the maps $d_\pm(s) =
(\pm s, s)$. Using the isotopy induced by the embedding of $\Sigma
D_i^\pm$, $\Sigma TD_i^\pm$ into $\NN'$ (we omit the details) we can
change the embedded spheres $S'_i, TS'_i$ to embedded spheres in a
general position given by the two formulas $S''_i = (D_i^+ \times 1)
\cup (S_i \times d_+I) \cup (D_i^- \times -1)$ and $TS''_i = (TD_i^-
\times -1) \cup (TS_i \times d_-I) \cup (TD_i^+ \times 1)$. It
follows that the associated intersection forms are $\lambda', \mu'$,
since the only intersections are the intersections in $\NN$.

The situation is depicted in Figure \ref{intersections}. The thicker
parts correspond to the spheres $S''_i$, $T S''_i$ in general
position. The only intersections are those in $\NN \times I$.
\end{proof}

\begin{prop} \label{prop2}
If $F^{(1)}(V,y)$ is not empty then the map $\ssigma_1 \colon
F^{(1)}(V,y) \rightarrow X_y$ is a homotopy equivalence.
\end{prop}

\begin{proof}
This is proved by showing that the map $\ssigma_1$ induces
isomorphisms on homotopy groups. Notice that we have
\[
\pi_k (F^{(1)}(V,y)) \cong \pi_k(X_y) \cong L_{k-n}(1).
\]
Because the surgery obstruction map $\theta$ is an isomorphism we
can realize any element of the group $L_{k-n}(1)$ by a $k$-simplex
in $F^{(1)}(V,y)$, that means by an equivariant splitting problem $f
\colon S(\VR) \times \Delta^{k+1} \rightarrow S(\VR) \times
\Delta^{k+1}$ along $S(V) \times \Delta^{k+1}$ satisfying certain
conditions on the boundary.

By Construction \ref{map3} the image $f' = \ssigma_1 (f)$ is a
$k$-simplex in $X_y$, that means an equivariant splitting problem $f
\colon S(\VR^2) \times \Delta^{k+1} \times I \rightarrow S(\VR^2)
\times \Delta^{k+1} \times I$ along $S(\VR) \times \Delta^{k+1}
\times I$ again with certain conditions on the boundary.

In order to show that the map $\ssigma_1$ induces an isomorphism on
homotopy groups it is enough to show that if the splitting
obstruction $\theta(f)$ is represented by a simple
$(-1)^m$-hermitian form $(H^+, \lambda_0,\mu_0)$, where $2m = k-n$,
then the splitting obstruction $\theta(f')$ is represented by the
same simple $(-1)^m$-hermitian form. But this is just a version of
the statement that was proved in a sequence of Lemmas \ref{lem1},
\ref{lem2}, \ref{lem3} parametrized over $\Delta^{k+1}$. The proofs
carry over to the parametrized case as well.
\end{proof}

\begin{thm} \label{thm2}
We have
\[
\pi_k(\TTheta F^{(1)}) \cong  L_k(1) \mathrm{\;\;for\;\;} k \in \Z,
\]
where the groups $L_k(1)$ are the $4$-periodic $L$-groups from
surgery theory associated to the trivial group.
\end{thm}

\begin{proof}
Notice that the proof of Proposition \ref{prop1} tells us in the
special case $y = \ast$ that modulo some low-dimensional deviations
the spectrum $\TTheta F^{(1)}$ is an $\Omega$-spectrum and we know
its homotopy groups. Namely
\[
\pi_k (\TTheta F^{(1)}) \cong \pi_{k+|k|+8}\big
(F^{(1)}(\R^{|k|+8})\big) \cong L_k(1)
\]
for all $k \in \Z$.
\end{proof}

%%%%%%%%%%%%%%%%%%%%%%%%%%%%%%%%%%%%%%%%%%%%%%%%%%

\

\

%%%%%%%%%%%%%%%%%%%%%%%%%%%%%%%%%%%%%%%%%%%%%%%%%%

\textbf{Acknowledgements.} The present paper grew out of my thesis
defended at the University of Aberdeen. I would like to thank my
supervisor Michael Weiss for conversations and encouragement during
the work on this project and my thesis referees Michael Crabb and
Andrew Ranicki for useful comments. My thanks are also to Mamuka
Jibladze for conversations and for the pictures in the paper.
Finally I would like to thank Max Planck Institut f\"ur Mathematik
in Bonn for the hospitality during final stages of the preparation
of this paper.

\end{document}